\documentclass[11pt,a4paper]{article}
\input amssymb.sty
\begin{document}
\newcommand{\Si}{\Sigma}
\newcommand{\tr}{{\rm tr}}
\newcommand{\ad}{{\rm ad}}
\newcommand{\Ad}{{\rm Ad}}
\newcommand{\ti}[1]{\tilde{#1}}
\newcommand{\om}{\omega}
\newcommand{\Om}{\Omega}
\newcommand{\de}{\delta}
\newcommand{\al}{\alpha}
\newcommand{\te}{\theta}
\newcommand{\vth}{\vartheta}
\newcommand{\be}{\beta}
\newcommand{\la}{\lambda}
\newcommand{\La}{\Lambda}
\newcommand{\D}{\Delta}
\newcommand{\ve}{\varepsilon}
\newcommand{\ep}{\epsilon}
\newcommand{\vf}{\varphi}
\newcommand{\vfh}{\varphi^\hbar}
\newcommand{\vfe}{\varphi^\eta}
\newcommand{\fh}{\phi^\hbar}
\newcommand{\fe}{\phi^\eta}
\newcommand{\G}{\Gamma}
\newcommand{\ka}{\kappa}
\newcommand{\ip}{\hat{\upsilon}}
\newcommand{\Ip}{\hat{\Upsilon}}
\newcommand{\ga}{\gamma}
\newcommand{\ze}{\zeta}
\newcommand{\si}{\sigma}
\def\bfa{{\bf a}}
\def\bfb{{\bf b}}
\def\bfc{{\bf c}}
\def\bfd{{\bf d}}
\def\bfe{{\bf e}}
\def\bff{{\bf f}}
\def\bfm{{\bf m}}
\def\bfn{{\bf n}}
\def\bfp{{\bf p}}
\def\bfu{{\bf u}}
\def\bfv{{\bf v}}
\def\bft{{\bf t}}
\def\bfx{{\bf x}}
\def\bfg{{\bf g}}
\def\bfC{{\bf C}}
\def\bfS{{\bf S}}
\def\bfJ{{\bf J}}
\def\bfr{{\bf r}}
\def\bfnu{{\bf \nu}}
\def\bfsi{{\bf \sigma}}
\def\bfU{{\bf U}}

\def\hS{{\hat{S}}}

\newcommand{\li}{\lim_{n\rightarrow \infty}}
\def\mapright#1{\smash{
\mathop{\longrightarrow}\limits^{#1}}}

\newcommand{\mat}[4]{\left(\begin{array}{cc}{#1}&{#2}\\{#3}&{#4}
\end{array}\right)}
\newcommand{\thmat}[9]{\left(
\begin{array}{ccc}{#1}&{#2}&{#3}\\{#4}&{#5}&{#6}\\
{#7}&{#8}&{#9}
\end{array}\right)}
\newcommand{\beq}[1]{\begin{equation}\label{#1}}
\newcommand{\eq}{\end{equation}}
\newcommand{\beqn}[1]{\begin{eqnarray}\label{#1}}
\newcommand{\eqn}{\end{eqnarray}}
\newcommand{\p}{\partial}
\newcommand{\di}{{\rm diag}}
\newcommand{\oh}{\frac{1}{2}}
\newcommand{\su}{{\bf su_2}}
\newcommand{\uo}{{\bf u_1}}
\newcommand{\SL}{{\rm SL}(2,{\mathbb C})}
\newcommand{\GLN}{{\rm GL}(N,{\mathbb C})}
\def\sln{{\rm sl}(N, {\mathbb C})}
\def\sl2{{\rm sl}(2, {\mathbb C})}
\def\SLN{{\rm SL}(N, {\mathbb C})}
\def\SLT{{\rm SL}(2, {\mathbb C})}
\newcommand{\gln}{{\rm gl}(N, {\mathbb C})}
\newcommand{\PSL}{{\rm PSL}_2( {\mathbb Z})}
\def\f1#1{\frac{1}{#1}}
\def\lb{\lfloor}
\def\rb{\rfloor}
\def\sn{{\rm sn}}
\def\cn{{\rm cn}}
\def\dn{{\rm dn}}
\newcommand{\rar}{\rightarrow}
\newcommand{\upar}{\uparrow}
\newcommand{\sm}{\setminus}
\newcommand{\ms}{\mapsto}
\newcommand{\bp}{\bar{\partial}}
\newcommand{\bz}{\bar{z}}
\newcommand{\bw}{\bar{w}}
\newcommand{\bA}{\bar{A}}
\newcommand{\bL}{\bar{L}}
\newcommand{\btau}{\bar{\tau}}

\newcommand{\Sh}{\hat{S}}
\newcommand{\vtb}{\theta_{2}}
\newcommand{\vtc}{\theta_{3}}
\newcommand{\vtd}{\theta_{4}}

\def\mC{{\mathbb C}}
\def\mZ{{\mathbb Z}}
\def\mR{{\mathbb R}}
\def\mN{{\mathbb N}}

\def\frak{\mathfrak}
\def\gg{{\frak g}}
\def\gJ{{\frak J}}
\def\gS{{\frak S}}
\def\gL{{\frak L}}
\def\gG{{\frak G}}
\def\gk{{\frak k}}
\def\gK{{\frak K}}
\def\gl{{\frak l}}
\def\gh{{\frak h}}
\def\gH{{\frak H}}

\newcommand{\ran}{\rangle}
\newcommand{\lan}{\langle}
\def\f1#1{\frac{1}{#1}}
\def\lb{\lfloor}
\def\rb{\rfloor}
\newcommand{\slim}[2]{\sum\limits_{#1}^{#2}}

\newcommand{\sect}[1]{\setcounter{equation}{0}\section{#1}}
\renewcommand{\theequation}{\thesection.\arabic{equation}}
\newtheorem{predl}{Proposition}[section]
\newtheorem{defi}{Definition}[section]
\newtheorem{rem}{Remark}[section]
\newtheorem{cor}{Corollary}[section]
\newtheorem{lem}{Lemma}[section]
\newtheorem{theor}{Theorem}[section]

\vspace{0.3in}
\begin{flushright}
 ITEP-TH-20/05\\
\end{flushright}
\vspace{10mm}
\begin{center}
{\Large{\bf Painlev{\'e} VI,
Rigid Tops and Reflection Equation}
}\\
\vspace{5mm}
A.M.Levin\\
{\sf Max Planck Institute of Mathematics, Bonn, Germany,}\\
{\sf Institute of Oceanology, Moscow, Russia,} \\
{\em e-mail alevin@wave.sio.rssi.ru}\\
M.A.Olshanetsky
\\
{\sf Max Planck Institute of Mathematics, Bonn, Germany,}\\
{\sf Institute of Theoretical and Experimental Physics, Moscow, Russia,}\\
{\em e-mail olshanet@itep.ru}\\
A.V.Zotov \\
{\sf Institute of Theoretical and Experimental Physics, Moscow, Russia,}\\
{\em e-mail zotov@itep.ru}\\
\vspace{5mm}
\end{center}

\begin{abstract}
 We show that the Painlev{\'e} VI equation has an equivalent form of
the non-autonomous  Zhukovsky-Volterra
gyrostat. This system is a generalization of the Euler top in $\mC^3$
and include the additional constant gyrostat momentum.
The quantization of its autonomous version is achieved
by the reflection equation.
The corresponding quadratic algebra generalizes the Sklyanin algebra.
 As by product we define integrable XYZ
spin chain on a finite lattice with new boundary conditions.
\end{abstract}

\today

\tableofcontents

\section{Introduction}
\setcounter{equation}{0}

In this paper we discuss a few issues related to isomonodromy problems on
elliptic curves, integrable systems with the spectral parameter on the
same curves and the XYZ spin-chain on a finite lattice.
Our main object is the  Painlev\'{e} VI
equation (PVI). It is a second order ODE depending on four free parameters
$(\al,\be,\ga,\de)$
$$
\frac{d^2X}{dt^2}=\frac{1}{2}\left(\frac{1}{X}+\frac{1}{X-1}+
\frac{1}{X-t}\right)
\left(\frac{dX}{dt}\right)^2-
\left(\frac{1}{t}+\frac{1}{t-1}+\frac{1}{X-t}\right)\frac{dX}{dt}+
$$
\beq{I.1}
+\frac{X(X-1)(X-t)}{t^2(t-1)^2}\left(\al+\be\frac{t}{X^2}+
\ga\frac{t-1}{(X-1)^2}
+\de\frac{t(t-1)}{(X-t)^2}\right).
\eq
PVI was discovered by B.Gambier \cite{Gam} in 1906.
He accomplished the Painlev\'{e}
classification program
of the second order differential equations whose solutions have not
movable critical points. This equation and its degenerations $PV-PI$ have
a lot of applications in Theoretical and Mathematical Physics.
(see, for example \cite{LW}).

We prove here that PVI can be write down in a very simple form as ODE
with a quadratic non-linearity. It is a
non-autonomous version of the $\SL$ {\it Zhukovsky-Volterra gyrostat} (ZVG)
\cite{Zhuk,Vol}. The ZVG generalizes the standard Euler top in the space $\mC^3$
by adding an external constant rotator momentum.
 The ZVG equation describes the evolution
of the momentum vector $\vec{S}=(S_1,S_2,S_3)$ lying on a $\SL$ coadjoint
orbit. We consider Non-Autonomous Zhukovsky-Volterra gyrostat (NAZVG)
\beq{ur3}
\p_\tau\vec{S}=\vec{S}\times (\vec{J}(\tau)\cdot\vec{S})+\vec{S}\times
\vec{\nu}'\,
\eq
where $\vec{J}(\tau)\cdot\vec{S}=(J_1S_1,J_2S_2,J_3S_3)$.
Three additional constants $\vec{\nu}'=(\nu_1',\nu_2',\nu_3')$
form the gyrostat momentum vector, and the vector $\vec{J}=\{J_\al(\tau)\}$,
$(\al=1,2,3)$ is the inverse inertia vector depending on the "time" $\tau$
in the following way. Let $\Si_\tau=\mC/(\mZ+\tau\mZ)$ be the elliptic curve, and $\wp(x,\tau)$
is the Weierstrass function. Then $J_\al(\tau)=\wp(\om_\al,\tau)$,
where $\om_\al$ are the half-periods $(\oh,\frac{\tau}2,\frac{1+\tau}2)$.
The constants $\nu_\al$ along with the value of the Casimir
function $\sum_\al S_\al^2=(\nu_0')^2$ are expressed through the four
constants $(\al,\be,\ga,\de)$ of the PVI. If the gyrostat momentum vanishes,
 $\vec{\nu}'=0$ the ZVG is simplified to {\it the non-autonomous Euler Top} (NAET) equation
\beq{eu}
\p_\tau\vec{S}=\vec{S}\times (\vec{J}(\tau)\cdot\vec{S})\,.
\eq

To establish the connection between the PVI eq. (\ref{I.1}) and (\ref{ur3})
we start with the elliptic form of (\ref{I.1})
\cite{Pain,Ma}
\beq{ur2}
\frac{d^2u}{d\tau^2}=
-\sum\limits_{a=0}^3\nu_a^2\wp'(u+\om_a,\tau)\,,~~(\om_a=(0,\om_\al))\,.
\eq
The interrelation between two sets of the constants $\nu_a$ and $\nu_a'$ is explained
in Sections 5.3 and 5.4.
 We call this equation EPVI to distinguish it from
(\ref{I.1}).
By replacing $\tau$ on the external time $t$
we come to the BC$_1$
{\it Calogero-Inozemtsev system} (CI)
\beq{ur1}
\frac{d^2u}{dt^2}=-\sum\limits_{a=0}^3\nu_a^2\wp'(u+\om_a,\tau)\,.
\eq
This equation is the simplest case of the integrable BC$_N$ CI hierarchy
\cite{In}.
 When the constants
are equal, (\ref{ur1}) describes the two-body
elliptic Calogero-Moser (CM) system in the center of the mass frame
\beq{ur0}
\frac{d^2u}{dt^2}=-\nu^2\wp'(2u,\tau)\,.
\eq
In this case EPVI (\ref{ur2}) assumes the form
\beq{r0}
\frac{d^2u}{d\tau^2}=-\nu^2\p_u\wp(2u,\tau)\,.
\eq

 The genuine interrelations between integrable
and isomonodromy hierarchies arise on the level of the corresponding linear
systems. The linear equations of the isomonodromy problem
look like a quantization (the Whitham  quantization) of the linear
problem for the integrable hierarchy \cite{LO}.
While for the integrable systems the Lax matrices are
sections of the Higgs bundles (one forms on the spectral curve),
they become the holomorphic components of the flat connections for the
monodromy preserving equations.
Recently, the Lax representation for (\ref{ur1})  was proposed in Ref.\,\cite{Z}.
It allows us to construct the linear system for EPVI (\ref{ur2}).

In our previous work \cite{LOZ} we proposed a  transformation of
the linear system for the N-body CM system to the linear system for
the integrable $\SLN$ elliptic
Euler-Arnold top \cite{STSR}. It is the so-called {\it Hecke correspondence}
of the Higgs bundles. It is accomplished by a singular gauge
transform ({\it  the modification}) of the Lax equations.
This action on the dynamical variables
 is a symplectomorphisms. In the simplest case it provides a
change of variables from the two particle CM (\ref{ur0}) to the
$\sl2$ autonomous Euler Top (ET)
\beq{ur00}
\p_t\vec{S}=\vec{S}\times (\vec{J}\cdot\vec{S})\,,~~(\vec{J}\cdot\vec{S})=
(J_\al S_\al)\,.
\eq
It opens a way to define the Lax matrix $L^{ET}(z)$ of the ET
from $L^{CM}$.
 The both models contain a free constant which is
the coupling constant for CM and the value of the Casimir function
of the $\SL$ coadjoint orbit for ET.
It should be mentioned that the analogous transformation
was used in \cite{Has,SklTak,Va} for other purposes.

In the similar way we prove that
 the autonomous $\SL$ ZVG
\beq{ur30}
\p_t\vec{S}=\vec{S}\times \hat{J}\vec{S}+\vec{S}\times
\vec{\nu}'\,.
\eq
  is derived from the BC$_1$ CI eq. (\ref{ur1}) via the same modification.
In this way we obtain the Lax matrix $L^{ZVG}(z)$ from the  Lax matrix
of the CI system.
The change of variables is given explicitly. For ${\rm SO}(3)$ ZVG
 the Lax pair with the rational spectral parameter was constructed
in Ref.\,\cite{Fedorov} and discussed in Ref.\,\cite{BoMa}.
In our case the Lax pair depends on the elliptic spectral parameter.

 In fact, the modification can be applied to the isomonodromy problem
\cite{AL,Ob}. It
acts on connections and, in particular, transforms
 (\ref{r0}) to  NAET (\ref{eu})
and the generic EPVI (\ref{ur2}) to the NAZVG
(\ref{ur3}). We present the explicit dependence $\vec{S}(u,\p_\tau u)$.
In this way we establish the equivalence between PVI (\ref{I.1}) and NAZVG
(\ref{ur3}).

There exists  another way to define the Lax matrix $L^{ZVG}(z)$ and thereby to derive ZVG
(\ref{ur30}) that will be used in this paper.
The starting point is a special Elliptic Garnier-Gaudin system (EGG).
EGG is an example of the Hitchin systems \cite{LOZ}.
It is derived from the rank two quasi-parabolic Higgs bundle
\cite{Si} over $\Si_\tau$ of degree one. We assume that the
Higgs field has simple poles at the half-periods $\om_a$, $(a=0,\ldots,3)$.
The invariant part of the Higgs field with respect to the
involution $z\to -z$ leads to $L^{ZVG}(z)$.

Finally, we can derive ZVG starting with the quantum reflection equation.
Consider first the quantization of $\SLT$ ET. It can be performed by
the quantum exchange relations with the Baxter $R$-matrix
\cite{Bax}:
\beq{ur4}
R(z,w)\hat{L}_1^{ET}(z)\hat{L}_2^{ET}(w)=
\hat{L}_2^{ET}(w)\hat{L}_1^{ET}(z)R(z,w)\,,
\eq
where $\hat{L}^{ET}(z)$ is the quantum Lax matrix of the elliptic top.
This equation is equivalent to the Sklyanin algebra \cite{Skl1}.
We show that the quantum Lax matrix for ZVG satisfies the reflection equation
introduced in \cite{Skl2}:
\beq{ur5}
R^-(z,w)\hat{L}_1^{ZVG}(z)R^+(z,w)\hat{L}_2^{ZVG}(w)=
\hat{L}_2^{ZVG}(w)R^+(z,w)\hat{L}_1^{ZVG}(z)R^-(z,w)\,.
\eq
The corresponding quadratic algebra generalizes the Sklyanin algebra.
In the classical limit (\ref{ur5}) yields two Poisson structures for ZVG:
\beq{ur6}
\{L_1(z),L_2(w)\}_2=\frac{1}{2}[L_1(z)L_2(w),r^-(z,w)]
-\frac{1}{2}L_1(z)r^+(z,w)L_2(w)+\frac{1}{2}L_2(w)r^+L_1(z)\,,
\eq
$$
\{L_1(z),L_2(w)\}_1=\frac{1}{2}[L_1(z)+L_2(w),r^-(z,w)]-
\frac{1}{2}[L_1(z)-L_2(w),r^+(z,w)]
$$
which are compatible as in the case of the Sklyanin algebra \cite{KLO}.
The first type of brackets generalizes the classical Sklyanin algebra while the second
is just Poisson-Lie brackets. The coefficient $\frac{1}{2}$ in (\ref{ur6})
comes from the statement that these brackets are derived from the standard
brackets
\beq{ur7}
\{L_1(z),L_2(w)\}_2=[L_1(z)L_2(w),r^-(z,w)]\,,
\eq
$$
\{L_1(z),L_2(w)\}_1=[L_1(z)+L_2(w),r^-(z,w)]
$$
by the Poisson reduction procedure for the constraints
$L(z)+L^{-1}(-z)\det L(-z)=0$  and $L(z)+L(-z)=0$
for the quadratic and linear brackets correspondingly.
This procedure however will not be discussed here.

In \cite{Skl2} the reflection equation was used to construct an
 integrable version of the XYZ spin-chain  on a finite lattice.
Following this recipe we obtain the XYZ spin-chain with the quantum ZVG on
the boundaries. Three additional constants here combine into the vector
of magnetic field. The classical Hamiltonian is presented in Proposition 8.3.
As far as we know the obtained model was not discussed earlier.

\bigskip
{\small {\bf Acknowledgments.}\\
 The work is supported by the grants NSh-1999-2003.2 of the scientific
schools, RFBR-03-02-17554  and CRDF RM1-2545. The work of A.Z. was
also partially supported by the grant MK-2059.2005.2.
We are grateful for hospitality
to  the Max Planck Institute of Mathematics at Bonn
where two of us (A.L. and M.O.) were staying
during preparation of this paper.
We would like to thank a referee  for  valuable remarks,
that allows us to improve the paper.}


\section{
Isomonodromic deformations and Elliptic Calogero-Moser
\\ System}
\setcounter{equation}{0}

We consider differential equations related to the $N$-body integrable
 elliptic Calogero-Moser system with spin
(CM) \cite{Ca,Mo,GH,Wo}.
 They are defined as monodromy preserving
equations of  some linear system on an elliptic curve and generalize
(\ref{ur0}) to $N$ dependent variables \cite{LO}. It
is a Hamiltonian non-autonomous system that describes dynamics of  $N$
particles with  internal degrees of freedom (spin)  in a
time-depending potential.
We call this system Non-autonomous Calogero-Moser system (NACM). In section
2.4 we consider the interrelations between autonomous equations,
corresponding to integrable
hierarchies and non-autonomous (monodromy preserving) equations.

\subsection{Phase space of NACM}

The phase space of NACM system is the same as of CM.
Let $\Si_\tau=\mC/\mZ^{(2)}_\tau$, $\mZ^{(2)}_\tau=\mZ\oplus\tau\mZ$,
$(\Im m\,\tau>0)$ be the elliptic curve.
The coordinates of the particles lie in $\Si_\tau$:
$$
\bfu=(u_1,\ldots,u_N)\,,~~u_j\in\Si_\tau
$$
with the constraint on the center of mass $~\sum u_j=0$.
Let $\mZ^{(2)}_\tau\ltimes W_N$ be the semi-direct product of the
two-dimensional lattice group and the permutation $W_N$.
The coordinate part of the phase space is the quotient
\beq{La}
\Lambda=(\mC^N/(\mZ^{(2)}_\tau\ltimes W_N))/\Si_\tau\,.
\eq
The last quotient respects the constraint on the center of mass.
Let
$$
\bfv=(v_1,\ldots,v_N)\,,~~v_j\in\mC\,,~\sum v_j=0\,.
$$
be the momentum vector dual to $\bfu\,:$ $\,\{v_j,u_k\}=\de_{jk}$.
The pair $(\bfv,\bfu)$ describes
the "spinless" part of the phase space.

The additional phase variables, describing the internal degrees of freedom
of the particles, are the matrix elements of the $N$-order
matrix $\bfp$. More exactly, we consider $\bfp$ as an
element of the Lie coalgebra $\sln^*$.
The linear (Lie-Poisson) brackets on
$\sln^*$  for the matrix elements have the form
\beq{3.1}
\{p_{jk},p_{mn}\}=p_{jn}\de_{km}-p_{mk}\de_{jn}\,.
\end{equation}
Let ${\cal O}$ be a coadjoint orbit
\beq{orb}
{\cal O}=\{\bfp\in\sln^*~|~\bfp=\hbox{Ad}^*_h \bfp^0,~h\in\SLN\,,~\bfp^0\in
 \gH^*\}\,,
\end{equation}
where $\gH$ is the Cartan subalgebra of $\sln$.
The phase space ${\cal R}^{CM}=T^*\La\times{\cal O}$ is a symplectic
manifold with the symplectic form
\beq{symp}
\om=\langle d\bfv\wedge d\bfu\rangle- \langle\bfp^0
dhh^{-1}dhh^{-1}\rangle\,,
\eq
where the brackets stand for the trace. The form $\om$ is invariant
with respect to the action of the diagonal subgroup $D=\exp\gH$:
$h\to hh_1$, $h_1\in D$. Therefore, we can go further and pass to the
symplectic quotient
\beq{torb}
\ti{\cal O}={\cal O}//D\,.
\eq
 It implies the following constraints:\\
 i) the moment constraint $p_{jj}=0$,\\
ii) the gauge fixing, for example, as $p_{j,j+1}=p_{j+1,j}$.\\
{\it Example.} Let $\bfp^0=\nu\di(N-1,-1,\ldots,-1)$. Then $\dim{\cal O}=2N-2$.
It is the most degenerate non-trivial orbit. It leads to the spinless model,
since in this case $\dim\ti{\cal O}=0$.
We should represent $\bfp^0$ in the special form that takes into account the
moment constraint (i):
\beq{J}
\bfp^0=\bfJ^C=\nu
\left(\begin{array}{ccccc}
0&1&1&\cdots&1\\
1&0&1&\cdots&1\\
\vdots&\vdots&\ddots&\ddots&\vdots\\
1&1&1&\cdots&1\\
1&1&1&\cdots&0
\end{array}\right)\,.
\eq
For $N=2$ these orbits are generic.

In this way we come to the phase space of the CM
\beq{b15}
{\mathcal R}^{CM}_{red} =\{T^*(\La)\times\ti{\cal O}\}\,,
\end{equation}
Note that
\beq{od}
\dim({\mathcal R}^{CM}_{red})=2N-2+\dim{\cal O}-2\dim(D)=\dim{\cal O}\,.
\end{equation}


\subsection{Equations of motion and Painlev{\'e} VI}

The CM Hamiltonian has the form
\beq{4.1}
H^{CM,spin}=\oh\sum_{j=1}^Nv_j^2-\sum_{j>k}p_{jk}p_{kj}E_2(u_j-u_k;\tau)\,,
\end{equation}
where $E_2(x;\tau) =\wp(x;\tau)+2\eta_1(\tau)$ is the second
Eisenstein function (\ref{A.2}) and $\tau$ plays the role of time.
\footnote
{In what follows we replace the Weierstrass function $\wp(x;\tau)$
 used in Introduction by the Eisenstein function $E_2(x;\tau)$.
It does not affect the equations of motion.}
For the orbit, corresponding to (\ref{J}), the spinless Hamiltonian is
\beq{4.1a}
H^{CM}=\oh\sum_{j=1}^Nv_j^2-\nu^2\sum_{j>k}E_2(u_j-u_k;\tau)\,.
\end{equation}

For general non-autonomous Hamiltonian systems it is convenient
to work with the extended phase space by including the time.
Here we deal with the extended space
${\cal R}^{ext}=({\cal R^{CM}},{\cal T})$,
where ${\cal T}=\{\tau\in\mC\,|\,\Im m\tau>0\}$.
Equip it with the degenerate two-form
$$
\om^{ext}=\om
-\f1{\ka}dH^{CM,spin}(\bfv,\bfu,\tau)\wedge d\tau\,,
$$
 where $\ka\geq 0$ is the
so-called {\it classical level}. Note that $\om^{ext}$ is invariant with
respect the modular transformations $\PSL$ of ${\cal T}$ \cite{LO}.
It means that $\om^{ext}$ can be restricted on the moduli space\\
${\cal M}={\cal T}/ \PSL$. The vector field
$$
{\cal V}_\tau=\sum_{j,k,l}\left(\frac{\p H^{CM,spin}}{\p u_j}\p_{v_j}-
\frac{\p H^{CM,spin}}{\p v_j}\p_{u_j}+\sum_{mn}
\frac{\p H^{CM,spin}}{\p p_{mn}}(p_{ml}\de_{nk}-p_{kn}\de_{lm})\p_{p_{kl}}
\right) +\ka\p_\tau
$$
annihilate $\om^{ext}$ and define the equations of motion of NACM system
$$
\frac{df(\bfv,\bfu,\bfp,\tau)}{d\tau}={\cal
V}_{\tau}f(\bfv,\bfu,\bfp,\tau)\,.
$$
In particular,
\beq{5.1}
\ka\p_\tau u_j=v_j\,,
\end{equation}
\beq{6.1}
\ka\p_\tau v_n=-\sum_{j\neq n}p_{jk}p_{kj}\p_{u_n}E_2(u_j-u_n;\tau)\,,
\end{equation}
\beq{7.1}
\ka\p_\tau\bfp=2[\bfJ_{\bfu}(\tau)\cdot\bfp,\bfp]\,,
\end{equation}
where the operator $\bfJ_\bfu\cdot\bfp$ is defined by the diagonal action
\beq{7a.1}
\bfJ_{\bfu}(\tau)\cdot\bfp~:~p_{jk}\,\to\,E_2(u_j-u_k;\tau)p_{jk}\,,
\eq
i.e. each matrix element $p_{jk}$ is multiplied on $E_2(u_j-u_k;\tau)$.
 For $N=2$ we put $u_1=-u_2=u$, $v_1=-v_2=v=\ka\p_\tau u$  and
come to the two body NACM model
\beq{P6}
\p^2_\tau u=-\nu^2\p_uE_2(2u)\,.
\eq
It coincides with (\ref{r0}) since $\p_u\wp(u)=\p_uE_2(u)$.


\subsection{Lax representation}

The goal of this subsection is the Lax representation of (\ref{5.1}) --
(\ref{7.1}) \cite{LO}.

\subsubsection{Deformation of elliptic curve}
Let
$T^2=\{(x,y)\in\mR\,|\,x,y\in\mR/\mZ\}$ be a torus.
 Complex structure on $T^2$ is defined by the complex coordinate
$$
\Si_{\tau_0}=\{z=x+\tau_0y\,,\Im m\,\tau_0>0\}\,,~~\,,~~
\Si_{\tau_0}\sim\mC/(\mZ+\tau_0\mZ)\,,
$$
or by the operator $\p_{\bz}$ annihilating the one form $dz$.

Consider the deformation of the complex structure that preserves the point
$z=0$.
Let $\chi(z,\bz)$ be the characteristic function of a
neighborhood of $z=0$.
For  two neighborhoods ${\cal U}'\supset{\cal U}$
 of  $z=0$ define the smooth function
\beq{cf}
\chi(z,\bz)=\left\{
\begin{array}{cl}
1,&\mbox{$z\in{\cal U}$ }\\
0,&\mbox{$z\in\Si_{\tau_0}\setminus {\cal U}'.$}
\end{array}
\right.
\eq
Consider a new complex coordinate $w$
\beq{ct}
  w=z-\frac{\tau-\tau_0}{\rho}(\bz-z)(1-\chi(z,\bz))\,,
~~(\rho=\tau_0-\btau_0)\,.
\eq
The new coordinate defines the deformed elliptic curve
$\Si_{\tau}=\{w=x+\tau y\}$.
The partial derivatives with respect to the new coordinates $(w,\bw)$ are
\beq{pd}
\left\{
\begin{array}{l}
\p_w=\frac{\bar{\al}}{|\al|^2-|\be|^2}
\left(\p_z+\bar{\mu}\p_{\bz}\right)\,,\\
\p_{\bw}=\frac{{\al}}{|\al|^2-|\be|^2}
\left(\p_{\bar{z}}+\mu\p_z\right)\,,
\end{array}
\right.
\eq
where
$$
\al=1+\frac{\tau-\tau_0}{\rho}(1-\chi(z,\bz))\,,~~
\be=-\frac{\tau-\tau_0}{\rho}(1-\chi(z,\bz))\,,
$$
and  $\mu$ is the {\it Beltrami differential}
$$
\mu=-\frac{\be}{\al}=\frac{\tau-\tau_0}{\rho}
\left(
\frac{1-\chi(z,\bz)-(z-\bz)\bp\chi(z,\bz)}
{1+\frac{\tau-\tau_0}{\rho}(1-\chi(z,\bz)-(z-\bz)\p\chi(z,\bz))}
\right)\,.
$$
Note, that
$$
\mu=
\left\{
\begin{array}{ll}
0\,, &\mbox{$z\in{\cal U}$ }\\
\frac{\tau-\tau_0}{\rho}\,  ,&\mbox{$z\in\Si_{\tau_0}\setminus {\cal U}'.$}
\end{array}
\right.
$$
The partial derivatives form the basis dual to the basis of the one-forms $(dw,d\bw)$.
The common prefactors in the r.h.s. in (\ref{pd}) is irrelevant in our construction.
For this reason we replace  the partial derivatives
by the vector fields that only annihilate
the corrseponding one-forms $(dw,d\bar{w})$.
 We preserve the same notations for them:
\beq{po1}
\left\{
\begin{array}{l}
\p_w=\p_z+\bar{\mu}\p_{\bz}\,,\\
\p_{\bw}=\p_{\bz}+\mu\p_z\,,
\end{array}
\right.
\eq

The pair $(\tau-\tau_0,\bar{\tau}-\bar{\tau}_0)$ plays the role of local
complex coordinates in the moduli space ${\cal M}$ of elliptic curves
with a marked point in a neighborhood of $\tau_0$.

In what follows we replace the complex  coordinates $(w,\bw)$
 by the   independent coordinates
$$
\left\{
\begin{array}{l}
 w=z-\frac{\tau-\tau_0}{\rho}(\bz-z)(1-\chi(z,\bz))\,,\\
\ti{w}=\bz\,,\\
\end{array}
\right.
$$
and the corresponding independent vector fields $(\p_{w},\p_{\ti{w}})$,
annihilating $d{\ti{w}}$ and $dw$
\beq{po}
\left\{
\begin{array}{l}
\p_{w}=\p_z\,,\\
\p_{\ti{w}}=\p_{\bz}+\mu\p_z\,.
\end{array}
\right.
\eq
Note that $\p_{\ti{w}}=\p_{\bw}$, while $\p_w$ is independent on $\bar{\mu}$ now.
We pass from $(w,\bw)$ to the chiral coordinates $(w,\ti{w}=\bz)$ because
only the holomorphic dependence on the moduli ${\cal M}$ is essential
in our construction.
The holomorphic coordinate $(\tau-\tau_0)$ plays the role of time in
the Hamiltonian systems we consider here (see (\ref{ur3}) and (\ref{ur2})).


\subsubsection{Flat bundles of degree zero}

Let $V_N^0$ be a flat vector bundle of rank $N$ and degree $0$
over the deformed elliptic curve $\Si_\tau$.
Consider the connections
$$
\left\{
\begin{array}{c}
\ka\p_w+L^{(0)}(w,\ti{w},\tau)\,,\\
 \p_{\ti{w}}+\bL^{(0)}(w,\ti{w},\tau)\,,
\end{array}
\right.
$$
$$
L^{(0)}(w,\ti{w},\tau)\,,~\bL^{(0)}(w,\ti{w},\tau)\in C^\infty Map(\Si_\tau,\sln)\,.
$$
The flatness of the bundle $V_N^0$ means
\beq{flat}
\p_{\ti{w}} L^{(0)}-\ka\p_w\bL^{(0)}+[\bL^{(0)},L^{(0)}]=0\,.
\eq

By means of the gauge transformations $f(w,\ti{w})\in C^\infty{\rm
Map}\,(\Si_\tau\to\GLN)$
$$
\bL^{(0)}\to f^{-1}\p_{\ti{w}}f+f^{-1}\bL^{(0)}f
$$
the connections of generic bundles of degree zero can
be put in the following form:
\\
$1.~\bL^{(0)}=0$. Then from the flatness (\ref{flat}) we have
$$
\p_{\ti{w}} L^{(0)}(w,\ti{w})=0\,.
$$
$2.$ The connection of generic bundles of deg$(V_N^0)=0$ have the following quasi-periodicity:
$$
L^{(0)}(w+1)=L^{(0)}(w)\,,~~
L^{(0)}(w+\tau)=\bfe(\bfu )L^{(0)}(w)
\bfe(-\bfu)\,,
$$
where the diagonal elements of
$\bfe(\bfu)=\di(\exp(2\pi i u_1),\ldots,\exp(2\pi i u_N))$
define {\it the moduli of holomorphic bundles}. We identify
$\bfu$ with the coordinates of particles.
In fact, $u_j,\,j=1,\ldots,N$ belong to the dual to $\Si_\tau$ elliptic
curve ({\it the Jacobian}), isomorphic to $\Si_\tau$.\\
$3.$ We assume that $L^{(0)}$ has a simple pole at $w=0$ and
$$
  Res|_{w=0}\, L^{(0)}(w)=\bfp\,.
$$

The conditions $1,2,3$ fix $L^{(0)}$ up to a diagonal matrix $P$
\beq{Lax}
L^{(0)}=P+X\,,
\eq
$$
P=\di(v_1,\ldots,v_N)\,,
$$
$$
X=\{X_{jk}\}\,,~(j\neq k)\,,~~
X_{jk}=p_{jk}\phi(u_j-u_k,w)\,.
$$
The function $\phi$ is determined by (\ref{A.3}).
The quasi-periodicity of $L^{(0)}$ is provided by (\ref{A.14}).
The free parameters $\bfv=(v_1,\ldots,v_N)$ of $P$ can be identified with
the momenta.

The flatness of the bundle upon the gauge transform amounts the
consistency of the system
\beq{cons}
\left\{
\begin{array}{ll}
  i. & (\ka\p_w+L^{(0)}(w,\tau))\Psi=0\,,\\
  ii. & \p_{\ti{w}}\Psi=0\,.
\end{array}
\right.
\eq
To come to the monodromy preserving equation,
 we assume that the Baker-Akhiezer vector $\Psi$ satisfies an
additional equation. Let ${\cal Y}$ be a monodromy matrix of
the system (\ref{cons}) corresponding to homotopically non-trivial cycles
$\Psi\to\Psi{\cal Y}$. The equation
\beq{mp}
\begin{array}{ll}
  iii. & (\ka\p_\tau+M^{(0)}(w))\Psi=0
\end{array}
\eq
means that $\p_\tau{\cal Y}=0$,
and thereby the monodromy is independent on the complex structure of $T^2$.
The consistency of i. and iii. is the monodromy preserving
equation
\beq{mpe}
\p_\tau L^{(0)}-\p_{w} M^{(0)} -\f1{\ka}[L^{(0)},M^{(0)}]=0\,.
\eq
In contrast with the standard Lax equation it has additional term
$\p_wM^{(0)}$.

\begin{predl}
The equation (\ref{mpe}) is equivalent to the monodromy preserving equations
(\ref{5.1}), (\ref{6.1}), (\ref{7.1}) for $L^{(0)}$ (\ref{Lax}) and
$$
M^{(0)}=\{Y_{jk}\}\,,~(j\neq k)\,,~~\di\,M^{(0)}=0
$$
$$
Y_{jk}=p_{jk}f(u_j-u_k,w)\,,~~
f(u,w)=\p_u\phi(u,w)\,.
$$
\end{predl}
Proof is based on the  Calogero functional equation (\ref{ad2}) and
the heat equation (\ref{A.4b}).


\subsection{Isomonodromic deformations and integrable systems}

We can consider the isomonodromy preserving equations as a deformation
({\it Whitham quantization}) of integrable equations \cite{Ta}.
The level $\ka$ plays the role of the deformation parameter.
Here we investigate the particular example -
the integrable limit of the vector generalization of PVI (\ref{5.1}) --
(\ref{7.1}) \cite{LO}.

Introduce the independent time
$t$ as $\tau=\tau_0+\ka t$ for $\ka\to 0$ and some fixed $\tau_0$.
It means that $t$ plays the role of a local coordinate in an neighborhood
of the point $\tau_0$ in the moduli space ${\cal M}$ of elliptic curves.
It follows from (\ref{ct}) that the limiting curve is
$\Si_{\tau_0}=\{z,\bz\}$.
In this limit we come to the equations of motion of CM (\ref{5.1}) -- (\ref{7.1}):
$$
\p_t u_j=v_j\,,
$$
\beq{CM}
\p_t v_n=-\sum_{j\neq n}p_{jk}p_{kj}\p_{u_n}E_2(u_j-u_n;\tau_0)\,,
\eq
$$
\p_t\bfp=2[\bfJ_{\bfu}(\tau_0) \cdot\bfp,\bfp]\,.
$$
The linear problem for this system is obtained from the linear problem
for the Isomonodromy problem (\ref{cons}), (\ref{mp})
by the analog of the quasi-classical limit in Quantum Mechanics.
Represent the Baker-Akhiezer function in the WKB form
\beq{WKB}
\Psi=\Phi\exp\bigl(\frac{{\cal S}^{(0)}}{\ka}+{\cal S}^{(1)}\bigr)\,.
\end{equation}
Substitute (\ref{WKB})
 in the linear system (\ref{cons}), (\ref{mp}). If
$\p_{\bar z}{\cal S}^{(0)}=0$ and $\p _t{\cal S}^{(0)}=0$,
then the terms of order $\ka^{-1}$ vanish.
In the quasi-classical limit we put
$\p{\cal S}^{(0)}=\la$.
In the zero order approximation we come to the linear system for CM
$$
\left\{
\begin{array}{ll}
 i. & (\la+L^{(0)}(z,\tau_0))Y=0\,,\\
  ii. & \p_{\bz}Y=0\,,\\
iii. & (\p_t+M^{(0)}(z,\tau_0))Y=0\,.
\end{array}
\right.
$$
The consistency condition of this linear system
$$
\p_tL^{(0)}(z)-[L^{(0)}(z),M^{(0)}(z)]=0\,,
$$
is equivalent to the Calogero-Moser
equations (\ref{CM}) \cite{Kr}.

The Baker-Akhiezer function $Y$ takes the form
$$
Y=\Phi e^{ t\frac{\p}{\p \tau_0}{\cal S}^{(0)}}\,.
$$

The same quasi-classical limit can be applied for the monodromy preserving
equations that will be considered in next Section.


\section{Isomonodromic deformations and Elliptic Top}
\setcounter{equation}{0}

\subsection{Euler-Arnold top on $\SLN$.}

Let $\bfS\in\sln^*$. Expand it in the basis (\ref{B.11})
$\bfS=\sum_{\al\in\ti{\mZ}^{(2)}_N}S_\al T_\al$.

According with (\ref{AA3b}) the Lie-Poisson brackets on $\sln^*$ assume the form
$$
\{S_\al,S_\be\}_1=\bfC(\al,\be) S_{\al+\be}\,.
$$
The Lie-Poisson
brackets are degenerated on $\gg^*=\sln^*$ and their symplectic leaves are
coadjoint orbits of $\SLN$. To descend to a particular coadjoint orbit
${\cal O}$ one should fix the values of the Casimirs for the linear bracket.
The phase space is a coadjoint orbit
\beq{cao}
{\cal R}^{ET}=\{\bfS\in\gg^*\,|\,\bfS=
gS_0 g^{-1}\,,~g\in\SLN\,,~S_0\in\gg^*\}\,.
\eq

The Hamiltonian of the Euler-Arnold top has a special form.
It is a quadratic functional on $\gg^*$
$$
H=-\oh\tr( \bfS,\bfJ(\bfS))\,,~~\bfS\in\gg^*\,,
$$
where $\bfJ$ is an invertible symmetric operator
$$
\bfJ\,:\,\gg^*\to\gg\,.
$$
It is called {\it the inverse inertia tensor}.
The equations of motion assume the form
\beq{EA}
\p_t\bfS=\{\bfJ(\bfS),\bfS\}_1=[\bfJ(\bfS),\bfS]\,.
\eq


\subsection{Non-autonomous Elliptic top (NAET).}

Consider a special form of the inverse inertia tensor.
Let
$$
E_{2}(\al)=E_2\bigl(\frac{\al_1+\al_2\tau}{N}|\tau\bigr)\,,
~~\al=(\al_1,\al_2)\in\ti{\mZ}^{(2)}_N\,,
$$
where $\ti{\mZ}^{(2)}_N$ is defined by (\ref{B.10}).
Define as above the diagonal action
$$
\bfJ(\tau,\bfS)~:~S_\al\to E_{2}(\al) S_\al\,,~~
\bfJ(\tau,\bfS)=\bfJ(\tau)\cdot\bfS\,.
$$
Then the Hamiltonian of ET assumes the form
$$
H^{NAET}(\bfS,\tau)=\frac{\pi\imath}{N^2}\tr(\bfS\bfJ(\tau)\cdot\bfS)
= -\oh\sum_{\ga\in\ti{\mZ}^{(2)}_N} S_\ga E_{2}(\ga) S_{-\ga}\,.
$$
 The equation of motion of NAET
is similar to (\ref{EA})
\beq{nat}
\ka\p_\tau\bfS=[\bfJ(\tau)\cdot\bfS,\bfS]\,,~~
(\,\ka\p_\tau S_\al=
\sum_{\ga\in\ti{\mZ}^{(2)}_N}S_\ga S_{\al-\ga}E_2(\ga)\bfC(\al,\ga)\,)\,.
\eq

As in the previous case one can consider the limit $\ka\to 0$,
$(\tau\to\tau_0)$
to the integrable elliptic top
\beq{iett}
\p_t\bfS=[\bfJ(\tau_0)\cdot\bfS,\bfS]\,.
\eq


\subsection{Lax representation}
\subsubsection{Flat bundles of degree one}

Let $V_N^1$ be a flat bundle over
the deformed elliptic curve $\Si_\tau$ of  rank $N$ and degree $1$
with the connections
\beq{lins}
\left\{
\begin{array}{l}
\ka\p_w+L^{(1)}(w,\ti{w},\tau)\,,\\
\p_{\ti{w}}+\bL^{(1)}(w,\ti{w},\tau)\,,
\end{array}
\right.
\eq
where $L^{(1)}(w,\ti{w},\tau)$, $\bL^{(1)}(w,\ti{w},\tau)$ are meromorphic maps of
$\Si_\tau$ in $\sln$.

For generic flat bundles of degree one the connections can be chosen in the form\\
$1.~\bL^{(0)}=0$. From the flatness  one has
$$
\p_{\ti{w}} L^{(1)}=0\,.
$$
$2.$ The Lax matrices satisfy the quasi-periodicity
conditions
$$
\begin{array}{c}
 L^{(1)}(w+1)=QL^{(1)}(w)Q^{-1}\,,  \\
 L^{(1)}(w+\tau)=\tilde\La L^{(1)}(w)\tilde\La^{-1}+\frac{2\pi i\ka}N\,,
\end{array}
$$
$$
\tilde\Lambda(w,\tau)=-\bfe_N\bigl(-w-\frac{\tau}2\bigr)\Lambda
$$
for $Q,\La$ (\ref{q}), (\ref{la}).
It means that there are no moduli parameters for  $E_N^1$.\\
$3.$ $L^{(1)}$ has a simple pole at $w=0$ and all degrees of freedom come from the residue
$$
Res|_{w=0}\, L^{(1)}(w)=\bfS\,.
$$
The Lax matrix is fixed by these conditions:
\begin{lem}
The connection assumes the form
\beq{l1}
L^{(1)}(w)=-\frac{\ka}{N}\p_w\ln\vth(w;\tau)Id+
\sum_{\al\in\ti{\mZ}_N^{(2)}} S_\al\vf_\al(w)T_\al\,.
\eq
where $\vf(\al,w)$ is defined by (\ref{vf}), and $T_\al$ are the basis
elements (\ref{B.11}).
\end{lem}
Fixing the connections we come from (\ref{lins}) to the linear system
\beq{li}
\left\{\begin{array}{ll}
 i. & (\ka\p_w+L^{(1)}(w))\Psi=0\,, \\
 ii &  \p_{\ti{w}}\Psi=0\,.
\end{array}
\right.
\eq
As above, the independence of the monodromy of (\ref{li}) means that the
Baker-Akhiezer vector satisfies the additional linear equation
\beq{mono}
iii.~(\ka\p_\tau+M^{(1)})\Psi=0\,.
\eq
\begin{lem}
The equation of motion of the non-autonomous top (\ref{nat})
is the monodromy preserving equation for (\ref{li}) with
 the Lax representation
\beq{ltop}
\p_\tau L^{(1)}-\p_w M^{(1)}+\f1{\ka}[M^{(1)},L^{(1)}]=0\,,
\eq
where $L^{(1)}$ is defined by (\ref{l1}),
$$
M^{(1)}=-\frac{\ka}N\p_\tau\ln\vth(w;\tau) Id+
\sum_{\ga\in\ti{\mZ}_N^{(2)}} S_\ga f_\ga(w) T_\ga \,,
$$
and $f_\ga(w)$ is defined by (\ref{f}).
\end{lem}
The proof of the equivalence of (\ref{nat}) and (\ref{ltop}) is based on the
 addition formula (\ref{ad2}) and the heat equation (\ref{A.4b})
 as in the case of CM.

In the quasi-classical limit $\ka\to 0$ we come to the integrable top on
$\SLN$ (\ref{iett}).


\section{Symplectic Hecke correspondence}
\setcounter{equation}{0}

We construct here a map from the phase space of
CM (\ref{b15}) to the phase space of ET (\ref{cao})
\beq{xi}
\Xi^+\,:\,{\cal R}^{CM}\to{\cal R}^{ET}\,,~~((\bfv,\bfu,\bfp)
 \mapsto\bfS)\,,
\eq
such that $\Xi^+$ is the symplectic map
$$
\Xi^{+*}\om(\bfS)=\om(\bfv,\bfu,\bfp)\,.
$$

To construct it we
define the map of the sheaves of sections $\G(V^{(0)}_N)\to\G(V_N^{(1)})$
such that it is an isomorphism on the complement to
$w$ and it has one-dimensional cokernel at $w\in\Sigma_\tau$~:
$$
0\longrightarrow \G(V^{(0)}_N)\stackrel{\Xi^+}{\longrightarrow}
\G(V^{(1)}_N)\longrightarrow \mC|_w\longrightarrow 0\,.
$$
It is the so-called {\it upper modification} of the bundle $E^{(0)}_N$
at the point $w$.

On the complement to the point $w$ consider the map
$$
\G(V^{(1)}_N)\stackrel{~\Xi^-}\longleftarrow \G(V^{(0)}_N)\,,
$$
such that $\Xi^-\Xi^+=$Id. It defines {\it the lower modification}
at the point $w$.

In general case the modifications lead to the {\it Hecke correspondence}
between the moduli spaces of the holomorphic bundles of degree $k$ and $k+1$
\beq{sh}
{\cal M}^k\rightarrow{\cal M}^{k+1}\,.
\eq
The modifications can be lifted to the Higgs bundle. They act as singular
gauge transformations on the Lax matrices and provide the
symplectomorphisms between Hitchin systems ({\it symplectic
Hecke correspondence}).

 The modifications on the Higgs bundles can be applied for
the monodromy preserving equation as well.
The action of the upper modification on the Lax matrices has the form
\beq{mod}
L^{(1)}=\Xi^+\ka\p\Xi^{+-1}+ \Xi^+ L^{(0)}\Xi^{+-1}\,.
\eq
This form of transformation provides its symplectic action.

Consider in details the map (\ref{xi}).
 According with its definition the upper modification $\Xi^+(z)$ is
characterized by
the following properties:\\
$\bullet$ {\it Quasi-periodicity}:
\beq{8.12}
\Xi^+(z+1,\tau)= -Q\times \Xi^+(z,\tau)\,,
\end{equation}
\beq{8.13}
\Xi^+(z+\tau,\tau)=\tilde\Lambda(z,\tau)\times
\Xi^+(z,\tau)
\times{\rm diag}
({\bf e}(u_j))\,,~~\tilde\Lambda(z,\tau)=-\bfe_N\bigl(-z-\frac{\tau}2\bigr)\Lambda
\end{equation}
$\bullet\bullet$ $\Xi^+(z)$ has a simple pole at $z=0$.
 Let $\bfr_i=(r_{i,1},\ldots,r_{i,N})$ be an eigen-vector  of the matrix
$\bfp\in\ti{\cal O}$ (\ref{torb}), $\bfp\bfr_i=p^0_i\bfr_i$. Then
$Res(\Xi)_{z=0} \bfr_i=0$.\\
The former condition provides that the quasi-periods of the transformed
Lax matrix corresponds to the bundle of degree one.
The latter condition implies that $L^{(1)}$ has only a simple pole at
$z=0$. The residue at the pole is identified with $\bfS$.

First, we construct $(N\times N)$-
matrix $\tilde\Xi(z,\bfu;\tau)$ that satisfies (\ref{8.12}) and
(\ref{8.13}) but has a special one-dimensional kernel:
\beq{8.14}
\tilde\Xi_{ij}(z, \bfu;\tau) =
\theta{\left[\begin{array}{c}
\frac{i}N-\frac12\\
\frac{N}2
\end{array}
\right]}(z-Nu_j, N\tau ),
\end{equation}
where $\theta{\left[\begin{array}{c}
a\\
b\end{array}
\right]}(z, \tau )$ is the theta function with a characteristic
$$
\theta{\left[\begin{array}{c}
a\\
b
\end{array}
\right]}(z , \tau )
=\sum_{j\in \Bbb Z}
{\bf e}\left((j+a)^2\frac\tau2+(j+a)(z+b)\right)\,.
$$

It can be proved that the kernel of $\tilde\Xi$ at $z=0$ is generated
by the following column-vector~:
$$
\left\{(-1)^l \prod_{j<k;j,k\ne l}
\vartheta(u_k-u_j,\tau)\right\}, \quad l=1,2,\cdots,N\,.
$$

Then the matrix $\Xi(z,\bfu,\bfr_i)$, $(\bfr_i=(r_{1,i},\ldots,r_{N,i}))$
assumes the form
\beq{xi1}
\Xi^+(z,\bfu,\bfr_i)=\tilde{\Xi}(z)\times{\rm diag}
\left(\frac{(-1)^l}{r_{l,i}}
\prod_{j<k;j,k\ne l}
\vartheta(u_k-u_j,\tau)\right).
\end{equation}
The transformation (\ref{mod}) with $\Xi$ (\ref{xi1}) leads to the map
${\cal R}^{CM}\to {\cal R}^{ET}$ (\ref{xi}).

For the spinless CM, defined by the coupling
constant $\nu^2$, this transformation leads to the degenerate orbit of
the NAET.

Note that the equation for the spin variables of CM
(\ref{7.1}) reminds the equation of motion for the NAET
 with the coordinate-dependent operator $\bfJ_\bfu$ (\ref{7a.1}).
The only difference is the structure of the phase spaces
${\mathcal R}^{CM}$ (\ref{b15}) and ${\cal R}^{ET}$ (\ref{cao}). The
upper modification $\Xi^+$ carries out the pass from ${\mathcal R}^{CM}$ to
${\cal R}^{ET}$. It depends only on the part of variables on
${\mathcal R}^{CM}$, namely on $\bfu$ and $\bfp$ through the eigenvector
$\bfr_j$.
For the rank two bundles it is possible to write down the explicit
dependence $\bfS(u,v,\nu)$. We postpone this example to the general case
of PVI in Sect.6.


\section{Elliptic Garnier-Gaudin models}
\setcounter{equation}{0}
\subsection{General construction}

Here we present another way to construct the elliptic form of PVI
(\ref{ur2}) that does not appeal to the transformation (\ref{tr}) of
its rational form (\ref{I.1}).
The advantage of this approach is a simple derivation of the linear system and the
corresponding Lax operator.

The  CM and ET hierarchies are particular case of the following
construction. The main objects are the Lax matrices and the
corresponding integrable systems.
Let
$\D=\{z_a,\ a=1...m\}$ be the divisor of marked points on $\Si_\tau$ and
$V_N^k(\Si_\tau\backslash \D)$ is a holomorphic vector bundle of rank $N$
and degree $k$
over $\Si_\tau$ with quasi-parabolic structure at the marked points \cite{Si}.
This data allows us to define an elliptic analog of the Garnier-Gaudin system.

The moduli space ${\cal M}^{(0|k|N)}$   of holomorphic vector bundles of
degree $k$ and rank $N$
over $\Si_\tau$  is parameterized by  $d$ elements
$\bfu=(u_1,\ldots,u_{d})$, $u_j\in\Si_\tau$ and
$$
d=\dim({\cal M}^{(0|k|N)})=\hbox{g.c.d.}(N,k)\,.
$$
The bundles $V_N^k(\Si_\tau\backslash \D)$ and
$V_N^{k+N}(\Si_\tau\backslash \D)$ are isomorphic.
 The isomorphisms is provided by  twisting by a linear
bundle \cite{At}. Thus, in the definition of ${\cal M}^{(0|k|N)}$ $k$ is defined up $(mod\, N)$.
The transition matrices for a generic  holomorphic vector bundle can be chosen
as follows. For $d=\hbox{g.c.d.}(N,k)$ define the diagonal matrix
$\bfU$ containing $l=N/d$ blocks $(u_1,\ldots,u_d)$
$$
\bfU=\di\left((u_1,\ldots,u_d),\ldots,(u_1,\ldots,u_d)\right)\,.
$$
Then the transition matrices corresponding to the two basic non-contractible cycles
take the form
$$
g_1=Q\,,~~g_\tau=-\bfe_N(-kz-\oh k\tau)\exp (2\pi\imath{\bf U})\La^k\,.
$$
The quasi-parabolic structure means the fixing a flag
structure $Fl_a$ at $z=z_a$. One can
 act on the set of flags by the group $D_k=Centr(g_1,g_\tau)$.
 It is easy to see that $D_k$ is a subgroup of diagonal matrices
  $(D_k\subseteq D\subset\SLN)$ and $\dim D_k=d$.
The moduli space of the quasi-parabolic
bundles is defined as
$$
{\cal M}^{(m|k|N)}\sim{\cal M}^{(0|k|N)}\times(\sqcap_a Fl_a)/D_k\,.
$$
Then
$$
\dim{\cal M}^{(m|k|N)}=\hbox{g.c.d.}(N,k)+\sum_{a=1}^m\dim Fl_a-
\hbox{g.c.d.}(N,k)=\sum_{a=1}^m\dim Fl_a\,.
$$
In what follows we assume that all flag varieties are non-degenerate.
Thereby
$$
\dim{\cal M}^{(m|k|N)}=\oh mN(N-1)\,.
$$

The elliptic Garnier-Gaudin models $EGG(m|k|N)$ are particular
examples of Hitchin
systems \cite{H} related to holomorphic vector bundles of
degree $k$ and rank $N$. Particular cases of the elliptic Garnier-Gaudin models  were
considered in Ref.\,\cite{N} $(k=0)$ and in Ref.\,\cite{STSR,SklTak1} $(k=1)$.

The Lax matrix $L^{(m|k|N)}(z)dz$ is a meromorphic section of the bundle
${\rm End}V_N^k\otimes\Om^{(1,0)}(\Si_\tau\backslash\D)$
with simple poles at $\{z_a\}\in \D$
and fixed residues
\beq{rs}
Res\,L^{(m|k|N)}(z_a)=\bfS^a\in{\rm sl}^*(N, {\mathbb C})\,.
\eq
More exactly, we assume that $\bfS^a\in {\cal O}^a$ is an
element of a non-degenerate coadjoint orbit.
 The Lax matrix
$L^{(m|k|N)}(z)dz$  plays the role of the
quasi-parabolic Higgs field in the Hitchin construction
\footnote{In what follows we omit the differential $dz$ in the
notation of the Lax matrices and work with $L(z)$ instead of $L(z)dz$}.

The Higgs field  respects the lattice action $\mZ\oplus\tau\mZ$:
\beq{*}
L^{(m|k|N)}(z+1)=g_1^{-1}L^{(m|k|N)}(z)g_1\,,\ \ \ \ \ \
L^{(m|k|N)}(z+\tau)=g_\tau^{-1} L^{(m|k|N)}(z)g_\tau\,.
\eq

 It follows from (\ref{*}) that the diagonal part of $L^{(m|k|N)}(z)$ is
a double-periodic meromorphic function. It depends on the additional
variables $(v_1,\ldots,v_k)$. The conditions (\ref{rs}), (\ref{*}) determine
the Lax matrix $L^{(m|k|N)}(z)$. Below we consider particular appropriated
cases.

The whole phase space
$$
{\cal R}^{(m|k|N)}=T^*{\cal M}^{(0|k|N)}\times({\cal O}^1,\ldots,{\cal O}^m)=
\{\bfv,\bfu;\bfp^1,\ldots,\bfp^m\}
$$
 is equipped with the Poisson brackets
$$
\{v_j,u_k\}=\de_{jk}\,,~~
\{p^a_{il},p_{jk}^b\}=\de^{ab}(p^a_{ik}\de_{lj}-p^a_{jl}\de_{ik})\,.
$$

In the case when $\dim\,({\cal M}^{(0|k|N)})\neq 0$ one can go further and
consider the symplectic quotient with respect to the defined above group $D_k$
\beq{rps}
{\cal R}^{(m|k|N)}_{red}={\cal R}^{(m|k|N)}//D_k\,.
\eq
The dimension of the reduced space
$$
\dim{\cal R}^{(m|k|N)}_{red}=\sum_{a=1}^m\dim{\cal O}^a=mN(N-1)
$$
is independent on the degree of the bundles.
In fact, all the spaces ${\cal R}^{(m|k|N)}_{red}$ for\\
 $k=0,\ldots,N-1$ are
 symplectomorphic \cite{LOZ}. The symplectomorphisms is provided
by the symplectic Hecke
correspondence placed between the Higgs bundles of degree $k$ and $k+1$ (\ref{sh}).

The commuting Hamiltonians  of $EGG(m|k|N)$ can be extracted from the
expansion over the basis of the double-periodic functions $E_i(z-z_a)$
(\ref{A.1}), (\ref{A.2}),   (\ref{A.2a})
$$
\f1{j}\tr(L^{(m|k|N)}(z))^j=\left\{
\begin{array}{ll}
\sum_aH^a_{11}E_1(z-z_a)\,, & (j=1)\\
H^0_{0j}+\sum_a\sum^j_{i=1}H^a_{ij}E_i(z-z_a)\,, & (j=2,\ldots,N)\,,
\end{array}
\right.
$$
where
$\sum_aH^a_{1j}=0$.
Equivalently, the commuting Hamiltonians can be defined by the integral
representations \cite{H}
$$
H^a_{ij}=\f1{j}\int_{\Si_\tau}\tr
(L^{(m|k|N)}(z))^j\mu^a_{ij}\,.
$$
Here $\mu^a_{ij}=\mu^a_{ij}(z,\bz)$ are $(1-j,1)$ differentials on $\Si_\tau$.
They are chosen to be
dual to the basis of the elliptic functions
$$
\f1{j}\int_{\Si_\tau}\mu^a_{ij}E_l(z-z_b)=\de^{ab}_{il}\,.
$$
They have the following form \cite{LO,LO2}.
For a marked point $z_a$
define  a smooth function $\chi_a(z,\bz)$ as for $z=0$ (\ref{cf})
Then for $a\neq 0$, $i\neq 0$
$$
\mu^a_{ij}=\frac{j}{2\pi\imath}(z-z_a)^{i-1}\bp\chi_a(z,\bz)
$$
and
$$
\mu_{02}^0=\frac{1}{\tau-\bar{\tau}}
\bp(\bz-z)
(1-\sum_{a=1}^m\chi_a(z,\bz))\,.
$$
The  higher order differentials $\mu_{0j}^0$, $j>2$ are related
to the so-called W-structures on $\Si_\tau$ \cite{LO2} and will not be discussed here.

The commuting Hamiltonians  generate the hierarchy of the Lax equations
\beq{LE}
\p_{(a,i,j)}L^{(m|k|N)}(z))=
[L^{(m|k|N)}(z)),M^{(m|k|N)}_{(a,i,j)}(z)]\,,~~\p_{(a,i,j)}=\p_{t_{(a,i,j)}}\,.
\eq
Here the matrices $M^{(m|k|N)}_{(a,i,j)}$ have the same quasi-periodicity
properties as $L^{(m|k|N)}$ and satisfy the equation
$$
[M^{(m|k|N)}_{(a,i,j)},\bfU]-\bp M^{(m|k|N)}_{(a,i,j)}=(L^{(m|k|N)})^{j-1}\mu^a_{ij}
-\p_{(a,i,j)}\bfU\,.
$$
The Lax equations (\ref{LE}) as well as the form of the M-matrices follow
from the equations of motion for the Hitchin systems. Their derivation can be found,
for example, in \cite{O}.


\subsection{Involution of the Higgs bundles}

Let $\varsigma$ be an involution $(\varsigma^2=Id)$ acting on the
space of sections of the Higgs bundle. Then we have decomposition
$$
\G({\rm End}V_N^k)=\G^+({\rm End}V_N^k)\oplus\G^-({\rm End}V_N^k)\,.
$$
Let
$$
L^{(m|k|N)\,\pm}=\oh\bigl(L^{(m|k|N)}\pm\varsigma (L^{(m|k|N)})\bigr)\,.
$$
Assume that the involution preserves the hierarchy
$$
\varsigma\{M^{(m|k|N)}_{a,i,j}\}=\{M^{(m|k|N)}_{a,i,j}\}\,,
$$
and we choose the invariant subset
$\{M^{(m|k|N)+}_{a,i,j}\}$.
Then the constraint
\beq{cons11}
L^{(m|k|N)\,-}(t)\equiv 0
\eq
is consistent with the involutions (\ref{LE})
$$
\p_{a,i,j}L^{(m|k|N)\,+}(z))=[L^{(m|k|N)\,+}(z)),M^{(m|k|N)+}_{a,i,j}(z)]\,.
$$

\bigskip

Consider the case of the special divisor
\beq{D}
\Delta=(\om_0=0\,,{\bf \om}_\al\,,~\al=1,2,3)\,.
\eq
Then the transformation $z\to -z$ preserves $\Si_\tau\setminus \Delta$. It can be
accompanied with an involution of $\sln$ to generate $ \varsigma$.
In next subsections we consider
two examples of the rank two bundles with the involutions of such type.


\subsection{Degree zero bundles}
We have derived the space ${\cal M}^{(0|0|N)}$ in 2.3.2. Remind that
it is described by
$N$ parameters $\bfu=(u_1,\ldots,u_N)$, $u_i\in\Si_\tau$, such that
$u_1+...+u_N=0$ and
the multiplicators are:
\beq{t1}
g_1=Id_N,\ \ \ \ \ g_\tau=\bfe({\bf -u})=
\hbox{diag}(\bfe(-u_1),...,\bfe(-u_N))\,,~~(\bfe(u)=\exp 2\pi\imath u)\,.
\eq
These conditions fix the Lax matrix of $EGG(m|0|N)$ up to a diagonal
matrix
$\di\,\bfv$,\\ $\bfv=(v_1,\ldots,v_N)$. It follows from (\ref{*}) and (\ref{t1})
that the Lax matrix assumes the form
\beq{t2}
L^{(m|0|N)}_{ij}(z)=\delta_{ij}(v_i+
\sum\limits_{a=1}^mp^a_{ii}E_1(z-z_a))+
(1-\delta_{ij})
\sum\limits_{a=1}^mp^a_{ij}\phi(u_i-u_j,z-z_a)\,.
\eq
The particular important case corresponding to a single marked point
was considered in Section 2.

Now consider the special form of (\ref{t2}) with $\Delta$  (\ref{D})
and $N=2$. The Lax matrix (\ref{t2}) can be gauge transform to the form
\beq{tt2}
L^{(4|0|2)}(z)=(v+\slim{a=0}{3}p_3^aE_1(z+\om_a))\si_3+
\eq
$$
\slim{a=0}{3}(p_+^a\bfe(2u\p_\tau\om_a)\phi(2u,z+\om_a)\si_+ +
p_-^a\bfe(-2u\p_\tau\om_a)\phi(-2u,z+\om_a)\si_-\,.
$$
Define the involution as
\beq{inv}
\begin{array}{l}
\varsigma^{(0)}(L^{(4|0|2)}(z))= -\si_1L^{(4|0|2)}(-z)\si_1\,,\\
\varsigma^{(0)}(M^{(4|0|2)}(z))= \si_1M^{(4|0|2)}(-z)\si_1\,.
\end{array}
\footnote{The difference of the $\varsigma$ action on $L$ and $M$ is due to
$L\in\Om^{(1,0)}(\Si_\tau)$ while $M\in\Om^{(0)}(\Si_\tau)$.}
\eq
\begin{predl}[\cite{Z}]
The invariant Lax matrix $L^{4|0|2\,+}$ with respect to  $\varsigma^{(0)}$
assumes the form
$$
L^{(4|0|2)\,+}(z)=\tilde{L}^{CI}(z)\,,
$$
\beq{t23}
\tilde{L}^{CI}=\mat{v}{0}{0}{-v}+\sum\limits_{a=0}^3
\tilde{L}_a^{CI},\ \ \ \ \
\tilde{L}_a^{CI}
=\mat{0}{\tilde{\nu}_a
\vf_a(2u,z+\om_a)}{\tilde{\nu}_a
\vf_a(-2u,z+\om_a)}{0}\,,
\eq
where $\tilde{\nu}^2_a=\frac{1}{2}\tr(\bfp^a)^2$,
$~\vf_a(2u,z+\om_a)=\bfe(-2u\p_\tau\om_a)\phi(2u,z+\om_a)$.
\end{predl}
{\it Proof}.\\
The projection of the Lax matrix on its anti-invariant part
should vanish.
It follows from (\ref{A.300}) and (\ref{A.14}) that
$$
\bfe(2u\p_\tau\om_a)\phi(2u,z+\om_a)~\mapright{(z\to -z)}~
-\bfe(-2u\p_\tau\om_a)\phi(-2u,z+\om_a)\,.
$$
Since $\si_1\si_3\si_1=-\si_3$, $\si_1\si_+\si_1=\si_-$
the vanishing of $L^{(4|0|2)\,-}(z)$ implies
\beq{cinv}
p_3^a=0\,,~~p_+^a=p_-^a\,,~(a=0,\ldots,3)\,.
\eq
Thus one can put $p_+^a=p_-^a=\tilde{\nu}_a$ and the invariant part
$L^{(4|0|2)\,+}(z)$ coincides with (\ref{t23}). $\Box$

\begin{rem}
Remind that upon the involution the phase space (\ref{rps})
is the symplectic quotient
${\cal R}^{(4|0|2)}_{red}={\cal R}^{(4|0|2)}//D_2$.
The conditions (\ref{cinv}) imply the moment constraint $\sum_ap_3^a=0$
and the gauge fixing $\sum_ap_+^a=\sum_ap_-^a$ and thereby provide the pass to
the symplectic quotient.
\end{rem}

According with Ref.\,\cite{Z}
the new Hamiltonian can be read off from the decomposition
\beq{ham}
\frac{1}{4}\tr(\tilde{L}^{CI})^2=
H^{CI}+ \frac{1}{2}\sum_{a=0}^3\tilde{\nu}^2_aE_2(z-\om_a)\,,
\eq
where
$$
H^{CI}=\frac{1}{2}v^2-\frac{1}{2}\sum_{a=0}^3{\nu}^2_aE_2(u-\om_a)\,,
$$
and
\beq{nut}
\left(
\begin{array}{l}
\tilde{\nu}_0
\\
\tilde{\nu}_1
\\
\tilde{\nu}_2
\\
\tilde{\nu}_3
\end{array}
\right)
=\frac{1}{2}\left(
\begin{array}{l}
1\ \ \ \ 1\ \ \ \ \ 1\ \ \ \ \ 1
\\
1\ \ \ \ 1\ \ -1\ -1
\\
1\ -1\ \ \ \ \ 1\ -1
\\
1\ -1\ -1\ \ \ \ \ 1
\end{array}\right)
\left(
\begin{array}{l}
\nu_0
\\
\nu_1
\\
\nu_2
\\
\nu_3
\end{array}
\right)
\eq
The canonical bracket $\{v,u\}=1$ and the Hamiltonian
$H^{CI}$ yields the BC$_1$ CI system (\ref{ur1}).

\begin{predl}[\cite{Z}]
There exists the matrix $M^{4|0|2\,+}=\tilde{M}^{CI}$,
$$
\tilde{M}^{CI}=\sum\limits_{a=0}^3 \tilde{M}_a^{CI},\ \ \ \ \
\tilde{M}_\al^{CI}
=\mat{0}{\tilde{\nu}_\al\vf'_a(2u,z+\om_\al)}{\tilde{\nu}_a
\vf'_a(-2u,z+\om_\al)}{0},
$$
such that the Lax equation is
equivalent to the BC$_1$ CI equation (\ref{ur1}).
\end{predl}

\subsection{Degree one bundles}
The Lax matrices in this case was  considered in \cite{STSR}.
The multiplicators of a bundle $V_N^1$ are:
$$
g_1=-Q,\ \ \ \
g_{\tau}=-\bfe_N\left(-\frac{\tau}{2}-
z\right)\La\,,~~\bfe_N=\exp\frac{2\pi\imath}{N}
$$
where $Q$ and $\La$ are the matrices
(\ref{q}), (\ref{la}).
There are no moduli parameters in this case.
The phase space  is a direct sum of the coadjoint orbits
${\mathcal O}^1\times ...\times{\mathcal O}^m$ and
${\cal R}^{(m|1|N)}={\cal R}^{(m|1|N)}_{red}$.
In the basis $\{T_\al\},\ (\al\in\ti{\mZ}^{(2)}_N)$ (\ref{B.11}) the Poisson
brackets are
\beq{t6}
\{S^a_\al,S^b_\be\}_1=\bfC(\al,\be)
\delta^{ab}S^a_{\al+\be}\,.
\eq

The Lax matrix is completely fixed by the quasi-periodicity conditions
\beq{t5}
L^{(m|1|N)}(z)=\sum\limits_{a=0}^{m-1}\sum\limits_{\al\in\ti{\mZ}^{(2)}_N}
T_\al S^a_{\al}\vf_\al(z-z_a)\,,
\ \ \ \ \vf_\al(z)=\bfe(z\p_\tau\om_\al)\phi(z,\om_\al)\,,\ \
\om_\al=\frac{\al_1+\al_2\tau}{N}\,.
\eq

Consider the case $N=2$ and $\Delta$ (\ref{D})
\beq{t28}
L^{(4|1|2)}(z)=\sum\limits_{a=0}^3
\sum\limits_{\al=1}^3S^a_\al\vf_\al(z-\om_a)\si_\al\,.
\eq
Define the involution
\beq{inv1}
\begin{array}{c}
\varsigma^{(1)}(L^{4|1|2}(z))=-L^{4|1|2}(-z)\,,\ \
\varsigma^{(1)}(M^{4|1|2}(z))=M^{4|1|2}(-z)\,.
\end{array}
\eq
\begin{predl}
The invariant Lax matrix $L^{(4|1|2)\,+}$ with $\Delta$ (\ref{D})
assumes the form
$$
L^{(4|1|2)\,+}(z)=L^{ZVG}(z)\,.
$$
\beq{t29}
L^{ZVG}=\sum\limits_{\al=1}^3
\left(
S_\al\vf_\al(z)+\tilde{\nu}_\al\vf_\al(z-\om_\al)
\right)\sigma_\al=
\sum\limits_{\al=1}^3\left(
S_\al\vf_\al(z)+{\nu'}_\al\frac{1}{\vf_\al(z)}
\right)\sigma_\al\,,
\eq
where $S_\al=S_\al^0$, $\tilde{\nu}_\al=S_\al^\al=\hbox{const}$ and
\beq{nu}
\nu'_\al=-\tilde{\nu}_\al
\bfe(-\om_\al\p_\tau\om_\al)\left(\frac{\vth'(0)}{\vth(\om_\al)}\right)^2\,.
\eq
\end{predl}
{\it Proof.}\\
It follows from (\ref{inv1}) and (\ref{sc}) that
$$
L^{(4|1|2)\,-}=\oh(L^{(4|1|2)}(z)-L^{4|1|2}(-z))=\sum_{\al\neq a}
S^\al_\al\varphi_\al(z-\om_a)\si_\al\,.
$$
If $L^{(4|1|2)\,-}=0$, then   $S_\al^a=0$ for $\al\neq a$. Fixing
 the Casimir functions before the involution as $\sum_\al (S_\al^a)^2=\ti{\nu}^2_a$
  we conclude that $S_\al^\al=\ti{\nu}_\al$.
     On the other hand
$$
L^{(4|1|2)\,+}=\sum_{\al=1}^3
\left(
S^0_\al\varphi_\al(z)+S^\al_\al\varphi_\al(z-\om_\al)
\right)\si_\al\,.
$$
and we come to the statement. $\Box$

Again one can construct Hamiltonians from the Lax matrix $L^{ZVG}$ as in
(\ref{ham}). Then the Hamiltonian
\beq{ZVG}
H^{ZVG}=-\frac{\pi^2}{2}\tr\left( \bfS,(\bfJ(\tau)\cdot\bfS+\bfnu')\right)
\eq
with the (\ref{t6}) Lie-Poisson brackets leads to the ZVG equation
(\ref{ur30}).

It is worthwhile to mention that we deal with three types of constants $\nu$, $\ti{\nu}$
and $\nu'$. The latter two are expressed through the former by (\ref{nut}) and (\ref{nu}).

\begin{predl}
There exists the matrix $M^{(4|1|2)\,+}=M^{ZVG}$,
\beq{t32}
M^{ZVG}(z)=\sum\limits_\al\left[
-S_\al\varphi_\al(z)(E_1(z+\om_\al)-E_1(\om_\al))+
\tilde{\nu}_\al\varphi(z-\om_\al)E_1(z)
\right]\sigma_\al=
\eq
$$
=\sum\limits_\al
-S_\al\frac{\varphi_1(z)\varphi_2(z)\varphi_3(z)}{\varphi_\al(z)}\sigma_\al
+E_1(z)L^{ZVG}(z)\,.
$$
such that the Lax equation is
equivalent to the autonomous $\SL$ ZVG
(\ref{ur20}).
\end{predl}
The proof of this statement essentially is the same as Proposition 5.2.


\section{Non-autonomous systems}
\setcounter{equation}{0}

Consider the non-autonomous systems NAZVG (\ref{ur3}), and
EPVI (\ref{ur2})
\beq{ur20}
\p_\tau\vec{S}=\vec{S}\times \vec{J}(\tau,\vec{S})+\vec{S}\times
\vec{\nu}'\,,
\eq
\beq{ur309}
\frac{d^2u}{d\tau^2}=
-\sum\limits_{a=0}^3\nu_a^2\wp'(u+\om_a,\tau)\,.
\eq
We remind the relation between the last equation and PVI (\ref{I.1}),
established in \cite{Ma}. Assume that the parameters of the equations
are related as follows
$$
\nu_0=\al\,,~\nu_1=-\beta\,,~\nu_2=\ga\,,~\nu_3=\oh-\de\,.
$$
Then the substitution in (\ref{I.1})
\beq{tr}
(u,\tau)\rar
\left(X=\frac{E_2(u|\tau)-e_1}{e_2-e_1}\,,~
t=\frac{e_3-e_1}{e_2-e_1}\right)\,,~~e_a=E_2(\om_a)
\eq
brings it in the form (\ref{ur309}).

The  Lax matrices for the non-autonomous systems can be obtained
from the Lax matrices for integrable systems, constructed in previous
section,
 by replacing the coordinates
$z\to w$ (\ref{ct}).
In this way we come to the following two statements.
\begin{predl}[\cite{Z}]
The equation of motion of EPVI (\ref{ur309}) has the Lax form
\beq{10.6}
\p_\tau L^{CI}-\p_wM^{CI}+\f1{\ka}[M^{CI},L^{CI}]=0
\eq
with
\beq{10.4}
L^{CI}=P+X\,,~~P=\di(v,-v)\,,
\eq
$$
X_{12}(u,w)=\sum_a\ti{\nu}_a\varphi_a
(2u,w+\om_a)\,,~~X_{21}(u,w)=X_{12}(-u,w)\,,
$$
\beq{10.5}
M^{CI}=Y_{jk}\,,~~(j\neq k\,,\,j,k=1,2)
\eq
$$
Y_{12}(u,w)=\sum_a\ti{\nu}_a
\varphi_a'(2u,w+\om_a)\,,~~Y_{21}(u,w)=Y_{12}(-u,w)\,.
$$
\end{predl}
\begin{predl}
The equation of motion of NAZVG (\ref{ur20}) has the Lax form with
\beq{10.8}
L^{NAZVG}=
-\frac{\ka}{2}\p_w\ln\vth(w;\tau)\si_0+
\sum_\al (S_\al\vf_\al(w)+\nu_\al\vf_\al(w-\om_\al))\si_\al \,.
\eq
\beq{10.9}
M^{NAZVG}=-\frac{\ka}{2}\p_\tau\ln\vth(w;\tau)\si_0
+\sum\limits_\al
-S_\al\frac{\varphi_1(w)\varphi_2(w)\varphi_3(w)}{\varphi_\al(w)}\sigma_\al
+E_1(w)L^{NAZVG}(\ka=0)\,.
\eq
\end{predl}
The proofs of the both statements  repeat the autonomous cases.
The only new ingredient is the heat equation (\ref{A.4b}), that should be used.

Now we establish the interrelation between these systems.
It is provided by the same modification transform as above. Namely
$$
L^{NAZVG}(\bfS,w)=\Xi^+\ka\p(\Xi^{+})^{-1}+ \Xi^+ L^{(CI)}\Xi^{+-1}\,,
$$
where $\Xi^+$ is defined by (\ref{xi1}).

The upper modification leads to the following relations
between the coordinates of the phase spaces
\beq{w241}
\begin{array}{l}
S_{1}=-v\frac{\vtb(0)}{\vth'(0)}\frac{\vtb(2u)}{\vth(2u)}
-\frac{\ka}{2}\frac{\vtb(0)}{\vth'(0)}\frac{\vtb'(2u)}{\vth(2u)}+
\\
\tilde{\nu}_0\frac{\vtb^2(0)}{\vtc(0)\vtd(0)}
\frac{{\vtc(2u)\vtd(2u)}}{\vth^2(2u)}
+\tilde{\nu}_1\frac{\vtb^2(2u)}{\vth^2(2u)}+
\tilde{\nu}_2\frac{\vtb(0)}{\vtd(0)}
\frac{\vtb(2u)\vtd(2u)}{\vth^2(2u)}+\tilde{\nu}_3\frac{\vtb(0)}{\vtc(0)}
\frac{\vtb(2u)\vtc(2u)}{\vth^2(2u)}\,,
\\
iS_{2}=v\frac{\vtc(0)}{\vth'(0)}\frac{\vtc(2u)}{\vth(2u)}+
\frac{\ka}{2}\frac{\vtc(0)}{\vth'(0)}\frac{\vtc'(2u)}{\vth(2u)}-
\\
\tilde{\nu}_0\frac{\vtc^2(0)}{\vtb(0)\vtd(0)}
\frac{{\vtb(2u)\vtd(2u)}}{\vth^2(2u)}
-\tilde{\nu}_1\frac{\vtc(0)}{\vtb(0)}
\frac{\vtc(2u)\vtb(2u)}{\vth^2(2u)}-
\tilde{\nu}_2\frac{\vtc(0)}{\vtd(0)}
\frac{\vtc(2u)\vtd(2u)}{\vth^2(2u)}-
\tilde{\nu}_3\frac{\vtc^2(2u)}{\vth^2(2u)}\,,
\\
S_{3}=-v\frac{\vtd(0)}{\vth'(0)}\frac{\vtd(2u)}{\vth(2u)}
-\frac{\ka}{2}\frac{\vtd(0)}{\vth'(0)}\frac{\vtd'(2u)}{\vth(2u)}+
\\
\tilde{\nu}_0\frac{\vtd^2(0)}{\vtb(0)\vtc(0)}
\frac{{\vtb(2u)\vtc(2u)}}{\vth^2(2u)}+
\tilde{\nu}_1\frac{\vtd(0)}{\vtb(0)}
\frac{\vtb(2u)\vtd(2u)}{\vth^2(2u)}+\tilde{\nu}_2
\frac{\vtd^2(2u)}{\vth^2(2u)}+
\tilde{\nu}_3\frac{\vtd(0)}{\vtc(0)}
\frac{\vtd(2u)\vtc(2u)}{\vth^2(2u)}\,.
\end{array}
\eq

The proof is based on the direct usage of (\ref{kz1}), (\ref{kz2}) and
Riemann identities for theta functions \cite{Mam}.


\section{Quadratic brackets and NAET}

\subsection{r-matrix structure}

The classical $r$-matrix is the quasi-periodic map from
$\Si_\tau$ to End$(E_N^{(1)})\otimes$End$(E_N^{(1)})$ \cite{BD}
\beq{6.1a}
r(w)=\sum_\ga\vf_\ga(w)T_\ga\otimes T_{-\ga}\,,
\eq
where $\vf_\ga$ is defined by (\ref{vf}).
It satisfies the classical
Yang-Baxter equation
$$
[r^{(12)}(w-w'),r^{(13)}(w)]+
[r^{(12)}(w-w'),r^{(23)}(w')]
$$
$$
+[r^{(13)}(w),r^{(23)}(w')]=0\,.
$$
By means of the $r$-matrix one can define the linear brackets.
Let $L^{(1)}(z)$ be the Lax matrix (\ref{l1}) for $\ka=0$
$$
L^{(1)}(w)=\sum_\ga S_\al\vf_\al(w)T_\al\,.
$$
\begin{predl}[\cite{KLO}]
The Lie-Poisson brackets on $\sln$
$$
\{S_\al,S_\be\}_1=\bfC(\al,\be)S_{\al+\be}
$$
 in terms of the Lax operator $L^{(1)}(\bfS,z)$
 are equivalent to the
following relation for the Lax operator
$$
\{L_1^{(1)}(w),L_2^{(1)}(w')\}_1=
[r(w-w'),L^{(1)}(w)\otimes Id+
Id\otimes L^{(1)}(w')]
$$
$$
L^{(1)}_1=L^{(1)}\otimes Id\,,~~
L^{(1)}_2=Id\otimes L^{(1)}\,.
$$
\end{predl}
The proof is based on the Fay three-section formula (\ref{ad3}).

\subsection{Quadratic Poisson algebra}

In addition to the  $N^2-1$ variables
$\bfS=\{S_\al,\,\al\in\ti{\mZ}_N^{(2)}\}$ introduce a new
variable $S_0$
and the  $\GLN$-valued Lax operator
$$
\ti{L}=-S_0Id+L^{(1)}(\bfS,w)\,.
$$
It satisfies the classical exchange algebra:
\beq{cer}
\{\ti{L}_1(w),\ti{L}_2(w')\}_2=[r(w-w'),\ti{L}_1(w)\otimes \ti{L}_2(w')]\,.
\eq
These brackets are Poisson, since the Jacobi identity is provided by the
classical Yang-Baxter equation.
\begin{predl}
The quadratic Poisson algebra (\ref{cer}) in the coordinates $(S_0,\bfS)$
takes the form
\beq{s0}
\{S_\al,S_0\}_2=
\sum_{\ga\neq\al}S_{\al-\ga}S_\ga
\bigl(E_2(\frac{\ga_1+\ga_2\tau}{N})-E_2(\frac{\al_1-\ga_1+(\al_2-\ga_2\tau)}{N})
\bigr)
\bfC(\al,\ga)\,,
\eq
$$
\{S_\al,S_\be\}_2=S_0S_{\al+\be}
\bfC(\al,\be)\\
+\sum_{\ga\neq\al,-\be}S_{\al-\ga}S_{\be+\ga}
\bff(\al,\be,\ga)\bfC(\ga,\al-\be)\,,
$$
where
$$
\bff(\al,\be,\ga)=
E_1(\frac{\ga_1+\ga_2\tau}{N})+
E_1(\frac{\be_1+\ga_1-\al_1+(\be_2+\ga_2-\al_2)\tau}{N})-
$$
$$
E_1(\frac{\be_1+\ga_1+(\be_2+\ga_2)\tau}{N})+
E_1(\frac{\al_1-\ga_1+(\al_2-\ga_2)\tau}{N})\,.
$$
\end{predl}
It is the classical Sklyanin-Feigin-Odesski (SFO) algebra
\cite{Skl,FO}.
These brackets are extracted from (\ref{cer}) by means of
(\ref{ir}), (\ref{ir1}).

\bigskip

Two Poisson structures are called {\it compatible}
(or, form {\it a Poisson pair}) if their linear combinations are Poisson
structures as well. It turns out that
the linear and quadratic Poisson brackets are compatible, namely,
there exists the one-parametric family of the Poisson brackets
$$
\{\bfS,\bfS\}_\la=\{\bfS,\bfS\}_2+\la\{\bfS,\bfS\}_1\,.
$$
In the case of the ET integrable hierarchy  these compatible
brackets provide the hierarchy with the bihamiltonian  structure \cite{KLO}.
The hierarchies of the monodromy preserving equations are
more intricate \cite{LO2} and we do not consider here the hierarchy  of NAET.
However, we have the following manifestation of the bihamiltonian structure.
\begin{predl}
In terms of the quadratic brackets the equation of motion of NAET
(\ref{nat}) assumes the form
$$
\ka\p_\tau S_\al=\{S_0,S_\al\}_2\,.
$$
\end{predl}
The proof follows immediately from (\ref{s0}). We replace the linear
brackets on quadratic and simultaneously the quadratic Hamiltonian
$H^{NAET}$ on the linear $S_0$.


\section{Reflection equation and generalized Sklyanin algebra}
\setcounter{equation}{0}

In this Section we present another Hamiltonian form of NAZVG (\ref{ur3}).
It is based on the quadratic Poisson brackets. The quantization of these
brackets is described by quantum reflection equation.

\subsection{Quantum reflection equation}

Let  $R^-$ be the quantum vertex R-matrix, that arises in the
XYZ model. We introduce also the matrix $R^+$
\beq{60}
R^\pm(z,w)=
\sum\limits_{a=0}^3\varphi^\frac{\hbar}{2}_\al(z\pm w)
\sigma_\al\otimes\sigma_\al
\eq
where $\vf_a^\frac{\hbar}{2}$ is defined by (\ref{90}).
Define the quantum Lax operator
\beq{61}
\hat{L}(z)=\hS_0\phi^\hbar(z)\si_0+
\sum_\al(\hS_\al\vf^\hbar_\al(z)+
\tilde{\nu}_\al\vf^\hbar_\al(z-\om_\al))\si_\al\,.
\eq
\begin{predl}
The Lax operator satisfies the quantum reflection equation
\beq{62}
R^{-}(z,w)\hat{L}_1(z)R^{+}(z,w)\hat{L}_2(w)=
\hat{L}_2(w)R^{+}(z,w)\hat{L}_1(z)R^{-}(z,w)\,.
\eq
iff its components $S_a$ generate the associative
algebra with relations
\beq{63}
[\tilde{\nu}_\al,\tilde{\nu}_\be]=0\,,~~[\tilde{\nu}_\al,\hS_a]=0\,,
\eq
\beq{64}
i[\hS_0,\hS_\al]_+=[\hS_\be,\hS_\ga]\,,
\eq
\beq{65}
[\hS_\ga,\hS_0]=i\frac{K_\be-K_\al}{K_\ga}[\hS_\al,\hS_\be]_+
-2i\f1{K_\ga}(\tilde{\nu}_\al\rho_\al\hS_\be-
\tilde{\nu}_\be\rho_\be\hS_\al)\,,
\eq
where
$$
K_\al=E_1(\hbar+\om_\al)-E_1(\hbar)-E_1(\om_\al)\,,~~
\rho_\al=-\exp(-2\pi \imath\om_\al\p_\tau\om_\al)\phi(\om_\al+\hbar,-\om_\al)\,.
$$
\end{predl}
For the proof see Appendix D.

If all $\ti{\nu}_\al=0$ the algebra (\ref{63}) -- (\ref{65}) coincides
with the Sklyanin algebra. Therefore,
it is a three parametric
deformation of the Sklyanin algebra.

\bigskip
Two elements
$$
C_1=\hS_0^2+\sum_{\al}\hS_\al^2\,,
$$
$$
C_2=\sum_{\al}
\hS_\al^2K_\al(K_\al-K_\be-K_\ga)+2\ti{\nu}_\al\rho_\al K_\al\hS_\al
$$
belong to the center of the generalized Sklyanin algebra (\ref{63}),
(\ref{64}).
They are the coefficients of the expansion of the quantum determinant
$$
\det_\hbar (\hat{L})=\tr P(\hat{L}(z,\hbar)\otimes\hat{L}^+(z,-\hbar))\,,
$$
where
$$
P=\sigma_0\otimes\sigma_0+\sum_\al\sigma_\al\otimes\sigma_\al\,,
$$
$$
\hat{L}^+(z,\hbar)=
\hS_0\phi^\hbar(z)\si_0-
\sum_\al(\hS_\al\vf^\hbar_\al(z)+
\tilde{\nu}_\al\vf^\hbar_\al(z-\om_\al))\si_\al\,.
$$
over the basis of the elliptic functions.
However, we do not know how to derive the
expression for the quantum determinant directly from the reflection equation.


\subsection{Classical reflection equations}

Consider (\ref{62}) in the limit $\hbar\to 0$. The classical $r^\pm$-matrices
are defined from the expansion
$$
R^\pm(z,w)=\left(\frac{\hbar}{2}\right)^{-1}
\si_0\otimes\si_0+r^\pm(z,w)+O(\hbar)\,,
$$
$$
r^\pm(z,w)=\sum_\al\vf_\al(z\pm w)\si_\al\otimes\si_\al\,.
$$
For the Lax operator one has
$$
\hat{L}(z)=\hbar^{-1}\hS_0\si_0+
\sum_\al(\hS_\al\vf_\al(z)+\tilde{\nu}_\al\vf_\al(z-\om_\al))\si_\al
+O(\hbar)\,.
$$
Define the classical variables
$$
\hS_\al\to S_\al\,,~~\hS_0\to \hbar S_0\,,
$$
and the corresponding classical Lax operator
$$
\ti{L}(z)=S_0\si_0+
\sum_{\al=1}^3(S_\al\vf_\al(z)+\nu_\al\vf_\al(z-\om_\al))\si_\al\,.
$$
Then, taking into account that
$[L_1,L_2]=\hbar\{\ti{L}_1,\ti{L}_2\}+O(\hbar^2)$
one finds the
classical reflection equation in the first order of $\hbar^{-1}$
\beq{66}
\{\tilde{L}_1(z),\tilde{L}_2(w)\}_2=
\frac{1}{2}[\tilde{L}_1(z)\tilde{L}_2(w),r^-(z,w)]+
\eq
$$
\frac{1}{2}\tilde{L}_2(w)r^+(z,w)\tilde{L}_1(z)-
\frac{1}{2}\tilde{L}_1(z)r^+(z,w)\tilde{L}_2(w)\,.
$$
by passing from the group-valued element $\ti{L}$. For the Lie-algebraic
element $L$
$$
 L(z)=\sum_{\al=1}^3
(S_\al\vf_\al(z)+\tilde{\nu}_\al\vf_\al(z-\om_\al))\si_\al
$$
we come to the linear brackets
\beq{67}
\{L_1(z),L_2(w)\}_1=-\frac{1}{2}
[r^-(z,w),L_1(z)+L_2(w)]
+\frac{1}{2}
[r^+(z,w),L_1(z)-L_2(w)]\,.
\eq
\begin{predl}
The classical reflection equations (\ref{66}) leads to the
quadratic Poisson structure on $\mC^4$
generalizing the classical Sklyanin algebra
\beq{68}
\{S_\al,S_\be\}_2=2i\ve_{\al\be\ga}S_0S_\ga\,,
\eq
\beq{69}
\{S_0,S_\al\}_2=i\ve_{\al\be\ga}S_\be S_\ga(E_2(\om_\be)-E_2(\om_\ga))+
2i\ve_{\al\be\ga}S_\be\nu'_\ga\,,
\eq
where $\nu'$ is defined as (\ref{nu}).
\end{predl}

{\it Proof.}\\
The only dissimilarity from the classical Sklyanin algebra is the linear term
in (\ref{69}). It comes from the last term in the right hand side of
(\ref{65}) in the limit $\hbar\to 0$. $\Box$

There are two Casimir elements of the Poisson algebra (\ref{68}), (\ref{69})
$$
c_1=\sum_\al S_\al^2\,,
$$
$$
c_2=S_0^2+\sum_\al(\wp(\om_\al)S_\al^2+2\nu_\al' S_\al)\,.
$$
They are the coefficients of the expansion of $\det(\ti{L}(z)$ over the
basis of elliptic functions.

\bigskip

It is easy to see that the brackets (\ref{68}) are compatible with the
linear sl$_2$ Lie-Poisson brackets
$$
\{S_\al,S_\be\}_1=2i\ve_{\al\be\ga}S_\ga\,.
$$
The equation of motion of NAZVG (\ref{ur3}) is written in terms of the linear
brackets:
$$
\p_\tau S_\al=\{H^{ZVG},S_\al\}_1\,.
$$
 The straightforward calculations shows that (\ref{ur3}) can be
written in the form
\beq{70}
\p_\tau S_\al=\{H_0,S_\al\}_2\,,~~H_0=S_0\,.
\eq
In this way for the generic form of PVI we have the same analog of the
bihamiltonian property as in the degenerate case (Proposition 5.4).


\subsection{Spin chain with boundaries}
Quantum reflection equation allows us to define XYZ model on a finite
lattice with  boundary conditions \cite{Skl1}. The Lax operator (\ref{61})
can be considered as a new solution of the reflection equation.

Consider a pair of matrices
$K^{\pm}(z)$ with Poisson brackets
\beq{bo1}
\{K_1^\pm,K_2^\pm\}=[K^\pm_1(z)K^\pm_2(w),r(z-w)]+K^\pm_2(w)r(z+w)K^\pm_1(z)-
K^\pm_1(z)r(z+w)K^\pm_2(w)
\eq
and $L^i(z),\ i=1...N$ with brackets
\beq{bo2}
\{L^i_1(z),L^j_2(w)\}=\delta^{ij}[r(z-w),L^i_1(z)L^j_2(w)]\,.
\eq
Then, according with \cite{Skl1}
\beq{bo3}
h(z)=tr\left[ K^+(z)L^N(z)...L^1(z)K^-(z)\left(L^1(-z)\right)^{-1}...
\left(L^N(-z)\right)^{-1}\right]
\eq
is the generating function of commuting Hamiltonians
\beq{bo4}
\{h(z), h(w)\}=0\,.
\eq
Choosing $L^{ZVG_\pm}(z)=L^{ZVG_\pm}(z,\pm\hbar)$ for
$K^\pm(z)$ we construct a spin chain with boundaries.
But in (\ref{66}) we have a factor $\frac{1}{2}$ which
comes from (\ref{60}) on the quantum level. Thus, we should put
the brackets on boundaries for $L^{ZVG_\pm}(z)$ to be two times more
than in (\ref{68}). In other words $R$-matrices for
$\hat{L}^{ZVG_\pm}(z)$ should depend on the same Planck constant
$\hbar$ as for all other $\hat{L}^i(z)$.

\begin{predl}
Spin chain involving $N$ internal vertices $L^i(z)$ with boundaries
$L^{ZVG_\pm}(z)$ is integrable if
\beq{bo5}
\begin{array}{l}
\{S^\pm_\al,S^\pm_\be\}=4\imath\varepsilon_{\al\be\ga}S^\pm_0S^\pm_\ga\,,
\\
\{S^\pm_0,S^\pm_\al\}=2\imath\varepsilon_{\al\be\ga}(S^\pm_\be
S^\pm_\ga\,,
(E_2(\om_\be)-E_2(\om_\ga))+S^\pm_\be\nu'_\ga)
\end{array}
\eq
and for $i,j=1...N$
\beq{bo6}
\begin{array}{c}
\{S^i_\al,S^j_\be\}=2\delta^{ij}\imath
\varepsilon_{\al\be\ga}S^i_0S^i_\ga\,,
\\
\{S^i_0,S^j_\al\}=\delta^{ij}\imath
\varepsilon_{\al\be\ga}S^i_\be S^i_\ga
(E_2(\om_\be)-E_2(\om_\ga))\,.
\end{array}
\eq
The nearest-neighbor interaction is described by the Hamiltonian:
\beq{bo7}
\begin{array}{c}
H=
\ln\left(S_0^-S_0^1+\sum\limits_\al S_\al^-S_\al^{1}(C-E_2(\om_\al))+
(\nu')_\al^- S_\al^1\right)+
\\
\ln\left(S_0^NS_0^++\sum\limits_\al
S_\al^NS_\al^{+}(C-E_2(\om_\al))+(\nu')_\al^+ S_\al^N\right)+
\\
\sum\limits_{i=1}^{N-1}\ln\left(
S_0^iS_0^{i+1}+\sum\limits_\al S_\al^iS_\al^{i+1}(C-E_2(\om_\al))\right)\,,
\end{array}
\eq
where $C$ is a constant, equal to the fraction of the
values of the Casimir functions for Sklyanin bracket (\ref{bo6}).
\end{predl}
The proof of (\ref{bo7}) is similar to those one given in \cite{FT2}
(Chapter III, $\S\, 5$).
Consider a special point $z_0\in\Si$, such that
$E_2(z_0)=\frac{S_0^iS_0^i+\sum\limits_{\al=1}^3S_\al^iS_\al^iE_2(\om_\al)}
{\sum\limits_{\al=1}^3S_\al^iS_\al^i}=C_i$ and assume that this relation is
independent on a point of the lattice
$C_i=C_j\ \forall i,j=1...N$. In this case
all $L^i(z_0)$ are degenerated and have the form
$L^i(z_0)=\al^i(z_0)\times\be^i(z_0)$ for some vectors
$\al^i(z_0)$ and covectors
$\be^i(z_0)$. To finish the proof we should notice that all $\det L^i(z)$
are Casimirs for (\ref{bo2}) and in our case
$L^i(z)L^i(-z)={1}\cdot\det L^i(z)\ \forall i.$

Some solutions of the reflection
equation corresponding to the spin chains with dynamical boundary conditions are constructed
in Ref.\,\cite{Kuznetsov, InKo}.


\section{Appendix}
\subsection{Appendix A. Elliptic functions.}
\setcounter{equation}{0}
\def\theequation{A.\arabic{equation}}

We assume that $q=\exp 2\pi i\tau$, where $\tau$ is the modular parameter
of the elliptic curve $E_\tau$.

The basic element is the theta  function:
\beq{A.1a}
\vth(z|\tau)=q^{\frac
{1}{8}}\sum_{n\in {\bf Z}}(-1)^n\bfe(\oh n(n+1)\tau+nz)=~~
(\bfe=\exp 2\pi\imath)
\eq
$$
q^{\frac{1}{8}}e^{-\frac{i\pi}{4}} (e^{i\pi z}-e^{-i\pi z})
\prod_{n=1}^\infty(1-q^n)(1-q^ne^{2i\pi z})(1-q^ne^{-2i\pi z})\,.
 $$
\bigskip

{\it The  Eisenstein functions}
\beq{A.1}
E_1(z|\tau)=\p_z\log\vth(z|\tau), ~~E_1(z|\tau)\sim\f1{z}-2\eta_1z,
\eq
where
\beq{A.6}
\eta_1(\tau)=\frac{3}{\pi^2}
\sum_{m=-\infty}^{\infty}\sum_{n=-\infty}^{\infty '}
\frac{1}{(m\tau+n)^2}=\frac{24}{2\pi i}\frac{\eta'(\tau)}{\eta(\tau)}\,,
\eq
where
$$
\eta(\tau)=q^{\frac{1}{24}}\prod_{n>0}(1-q^n)\,.
$$
is the Dedekind function.
\beq{A.2}
E_2(z|\tau)=-\p_zE_1(z|\tau)=
\p_z^2\log\vth(z|\tau),
~~E_2(z|\tau)\sim\f1{z^2}+2\eta_1\,.
\eq
The highest Eisenstein functions
\beq{A.2a}
E_j(z)=\frac{(-1)^j}{(j-1)!}\p^{(j-2)}E_2(z)\,,~~(j>2)\,.
\eq

{\it Relation to the Weierstrass functions}
\beq{a100}
\zeta(z,\tau)=E_1(z,\tau)+2\eta_1(\tau)z\,,
\eq
\beq{a101}
\wp(z,\tau)=E_2(z,\tau)-2\eta_1(\tau)\,.
\eq

The next important function is
\beq{A.3}
\phi(u,z)=
\frac
{\vth(u+z)\vth'(0)}
{\vth(u)\vth(z)}\,.
\eq
\beq{A.300}
\phi(u,z)=\phi(z,u)\,,~~\phi(-u,-z)=-\phi(u,z)\,.
\eq
It has a pole at $z=0$ and
\beq{A.3a}
\phi(u,z)=\frac{1}{z}+E_1(u)+\frac{z}{2}(E_1^2(u)-\wp(u))+\ldots\,.
\eq
Let
\beq{A3c}
f(u,z)=\p_u\phi(u,z)\,.
\eq
Then
\beq{A3b}
f(u,z)=\phi(u,z) (E_1(u+z)-E_1(u))\,.
\eq

{\it Heat equation}
\beq{A.4b}
\p_\tau\phi(u,w)-\f1{2\pi i}\p_u\p_w\phi(u,w)=0\,.
\eq

{\it Quasi-periodicity}

\beq{A.11}
\vth(z+1)=-\vth(z)\,,~~~\vth(z+\tau)=-q^{-\oh}e^{-2\pi iz}\vth(z)\,,
\eq
\beq{A.12}
E_1(z+1)=E_1(z)\,,~~~E_1(z+\tau)=E_1(z)-2\pi i\,,
\eq
\beq{A.13}
E_2(z+1)=E_2(z)\,,~~~E_2(z+\tau)=E_2(z)\,,
\eq
\beq{A.14}
\phi(u,z+1)=\phi(u,z)\,,~~~\phi(u,z+\tau)=e^{-2\pi \imath u}\phi(u,z)\,.
\eq
\beq{A.15}
f(u,z+1)=f(u,z)\,,~~~f(u,z+\tau)=e^{-2\pi \imath u}f(u,z)-2\pi\imath\phi(u,z)\,.
\eq

 {\it  The Fay three-section formula:}

\beq{ad3}
\phi(u_1,z_1)\phi(u_2,z_2)-\phi(u_1+u_2,z_1)\phi(u_2,z_2-z_1)-
\phi(u_1+u_2,z_2)\phi(u_1,z_1-z_2)=0\,.
\eq
Particular cases of this formula is
the  Calogero functional equation
\beq{ad2}
\phi(u,z)\p_v\phi(v,z)-\phi(v,z)\p_u\phi(u,z)=(E_2(v)-E_2(u))\phi(u+v,z)\,,
\eq

Another important relation is
\beq{ir}
\phi(v,z-w)\phi(u_1-v,z)\phi(u_2+v,w)
-\phi(u_1-u_2-v,z-w)\phi(u_2+v,z)\phi(u_1-v,w)=
\eq
$$
\phi(u_1,z)\phi(u_2,w)f(u_1,u_2,v)\,,
$$
where
\beq{ir1}
f(u_1,u_2,v)=\ze(v)-\ze(u_1-u_2-v)+\ze(u_1-v)-\ze(u_2+v)\,.
\eq
One can rewrite the last function as
\beq{ir3}
f(u_1,u_2,v)=-\frac{
\vth'(0)\vth(u_1)\vth(u_2)\vth(u_2-u_1+2v)
}{
\vth(u_1-v)\vth(u_2+v)\vth(u_2-u_1+v)\vth(v)
}\,.
\eq

\bigskip
{\it Theta functions with characteristics:}\\
For $a, b \in \Bbb Q$ by definition:
$$\theta{\left[\begin{array}{c}
a\\
b
\end{array}
\right]}(z , \tau )
=\sum_{j\in \Bbb Z}
{\bf e}\left((j+a)^2\frac\tau2+(j+a)(z+b)\right)\,.
$$
In particular, the function $\vth$ (\ref{A.1a}) is a
theta function with characteristics:
\beq{A.29}
\vartheta(x,\tau)=\theta\left[
\begin{array}{c}
1/2\\
1/2
\end{array}\right](x,\tau)\,.
\eq
Properties:
$$
\theta{\left[\begin{array}{c}
a\\
b
\end{array}
\right]}(z+1 , \tau )={\bf e}(a)
\theta{\left[\begin{array}{c}
a\\
b
\end{array}
\right]}(z  , \tau )\,,
$$
$$
\theta{\left[\begin{array}{c}
a\\
b
\end{array}
\right]}(z+a'\tau , \tau )
={\bf e}\left(-{a'}^2\frac\tau2 -a'(z+b)\right)
\theta{\left[\begin{array}{c}
a+a'\\
b
\end{array}
\right]}(z , \tau )\,,
$$
$$\theta{\left[\begin{array}{c}
a+j\\
b
\end{array}
\right]}(z , \tau )=
\theta{\left[\begin{array}{c}
a\\
b
\end{array}
\right]}(z , \tau ),\quad j\in \Bbb Z\,.
$$
The following notations are used:
$\theta\left[\begin{array}{l}a/2\\b/2\end{array}\right]=\theta_{ab}$
and $\vth=\theta_{11}$.

\subsection{Appendix B.  Lie algebra $\sln$ and elliptic functions}
\setcounter{equation}{0}
\def\theequation{B.\arabic{equation}}

Introduce the notation
$$
{\bf e}_N(z)=\exp (\frac{2\pi i}{N} z)
$$
 and two matrices
\beq{q}
Q=\di({\bf e}_N(1),\ldots,{\bf e}_N(m),\ldots,1)
\eq
\beq{la}
\La=
\left(\begin{array}{ccccc}
0&1&0&\cdots&0\\
0&0&1&\cdots&0\\
\vdots&\vdots&\ddots&\ddots&\vdots\\
0&0&0&\cdots&1\\
1&0&0&\cdots&0
\end{array}\right)\,.
\eq
Let
\beq{B.10}
\mZ^{(2)}_N=(\mZ/N\mZ\oplus\mZ/N\mZ)\,,~~\ti{\mZ}^{(2)}_N=
\mZ^{(2)}_N\setminus(0,0)
\eq
be the two-dimensional lattice of order $N^2$ and
$N^2-1$ correspondingly.
The matrices $Q^{a_1}\La^{a_2}$, $a=(a_1,a_2)\in\mZ^{(2)}_N$
generate a basis in the group $\GLN$, while $Q^{\al_1}\La^{\al_2}$,
$\al=(\al_1,\al_2)\in\ti{\mZ}^{(2)}_N$ generate a basis in the Lie algebra $\sln$.
Consider the projective representation of
$\mZ^{(2)}_N$ in $\GLN$
\beq{B.11}
a\to T_{a}=
\frac{N}{2\pi i}\bfe_N(\frac{a_1a_2}{2})Q^{a_1}\La^{a_2}\,,
\eq
\beq{AA3a}
T_aT_b=\frac{N}{2\pi i}\bfe_N(-\frac{a\times b}{2})T_{a+b}\,, ~~
(a\times b=a_1b_2-a_2b_1)
\eq
Here $\frac{N}{2\pi i}
\bfe_N(-\frac{a\times b}{2})$ is a non-trivial two-cocycle
in $H^2(\mZ^{(2)}_N,\mZ_{2N})$.
It follows from (\ref{AA3a}) that
\beq{AA3b}
[T_{\al},T_{\be}]=\bfC(\al,\be)T_{\al+\be}\,,
\eq
where
$\bfC(\al,\be)=\frac{N}{\pi}\sin\frac{\pi}{N}(\al\times \be)$ are
 the structure constants of $\sln$.

Introduce the following  constants on $\ti{\mZ}^{(2)}$:
\beq{AA50}
\vth(\ga)=\vth\bigl(\frac{\ga_1+\ga_2\tau}{N}\bigr)\,,
\eq
\beq{AA5}
E_1(\ga)=E_1\bigl(\frac{\ga_1+\ga_2\tau}{N}\bigr)\,,
~~E_2(\ga)=E_2\bigl(\frac{\ga_1+\ga_2\tau}{N}\bigr)\,,
\eq
and the quasi-periodic functions on $\Si_\tau$
\beq{phi}
\phi_\ga(z)=\phi(\frac{\ga_1+\ga_2\tau}{N},z)\,,
\eq
\beq{vf}
\vf_\ga(z)=\bfe_N(\ga_2z)\phi_\ga(z)\,,
\eq
\beq{f}
f_\ga(z)=
\bfe_N(\ga_2z)\p_u\phi(u,z)|_{u=\frac{\ga_1+\ga_2\tau}{N}}\,.
\eq
\beq{fe}
f_\ga(z)=\bfe_N(\ga_2z)\phi_\ga(z)(E_1(\frac{\ga_1+\ga_2\tau}N+z)-
E_1(\frac{\ga_1+\ga_2\tau}N))\,.
\eq

It follows from (\ref{A.3}) that
\beq{qp}
\vf_\ga(z+1)=\bfe_N(\ga_2)\vf_\ga(z)\,,~~
\vf_\ga(z+\tau)=\bfe_N(-\ga_1)\vf_\ga(z)\,.
\eq
\beq{qpf}
f_\ga(z+1)=\bfe_N(\ga_2)f_\ga(z)\,,~~
f_\ga(z+\tau)=\bfe_N(-\ga_1)f_\ga(z)-2\pi\imath\vf_\ga(z)\,.
\eq

\bigskip
{\it $\SLT$ case}\\
For   $\SLT$  instead of $T_\al$ we use the basis of sigma-matrices
\beq{100}
\si_0=Id\,,~~\si_1=i\pi T_{0,1}\,,~~\si_2=i\pi T_{1,1}\,,~~\si_3=-i\pi
T_{1,0}\,,
\eq
$$
\{\si_a\}=\{\si_0,\si_\al\}\,,(a=0,\al)\,,(\al=1,2,3)
$$
$$
\si_+=\frac{\si_1-\imath\si_2}2\,,~~\si_-=\frac{\si_1+\imath\si_2}2\,.
$$
The standard theta-functions with the characteristics are
\beq{101}
\te_{0,0}=\te_3\,,~~\te_{1,0}=
\te_2\,,~~\te_{0,1}=\te_4\,,~~\te_{1,1}=\te_1\,.
\eq

\begin{table}
\begin{center}
\begin{tabular}{|c||c|c|c|}
\hline
$\al$ & (1,0) & (0,1) & (1,1) \\
\hline
$\si_\al$ &
$\si_3$   &
$\si_1$ &
$\si_2$ \\
\hline
 half-periods & $\om_1=\oh$ & $\om_2=\frac{\tau}{2}$ &
$\om_3=\frac{1+\tau}{2}$ \\
\hline
$\varphi_\al (z)$    &
$\frac{\te_2(z)\te'_1(0)}{\te_2(0)\te_1(z)}$ &
$\frac{\te_4(z)\te'_1(0)}{\te_4(0)\te_1(z)}$ &
$\frac{\te_3(z)\te'_1(0)}{\te_3(0)\te_1(z)} $\\
\hline
\end{tabular}
\end{center}
\end{table}

\beq{efi02}
\vf_\al(z)\vf_\al(z-\om_\al)=-
\bfe_1(-\om_\al\p_\tau\om_\al)\left(\frac{\vth'(0)}{\vth(\om_\al)}\right)^2
\eq

\beq{sc}
\vf_a(-z-\om_b)=-\vf_a(z-\om_a)\de_{ab}+(1-\de_{ab})\vf_a(z-\om_b)\,.
\eq

{\it Formulae with doubled modular parameter:}
\beq{kz1}
\begin{array}{l}
\vtd(x,\tau)\vtc(y,\tau)+\vtd(y,\tau)\vtc(x,\tau)=
2\vtd(x+y,2\tau)\vtd(x-y,2\tau)
\\
\vtd(x,\tau)\vtc(y,\tau)-\vtd(y,\tau)\vtc(x,\tau)=
2\vth(x+y,2\tau)\vth(x-y,2\tau)
\\
\vtc(x,\tau)\vtc(y,\tau)+\vtd(y,\tau)\vtd(x,\tau)=
2\vtc(x+y,2\tau)\vtc(x-y,2\tau)
\\
\vtc(x,\tau)\vtc(y,\tau)-\vtd(y,\tau)\vtd(x,\tau)=
2\vtb(x+y,2\tau)\vtb(x-y,2\tau)
\end{array}
\eq

\beq{kz2}
\begin{array}{l}
2\vth(x,2\tau)\vtd(y,2\tau)=\vth(\frac{x+y}{2},\tau)\vtb(\frac{x-y}{2},\tau)
+\vtb(\frac{x+y}{2},\tau)\vth(\frac{x-y}{2},\tau)
\\
2\vtc(x,2\tau)\vtb(y,2\tau)=\vth(\frac{x+y}{2},\tau)\vth(\frac{x-y}{2},\tau)
+\vtb(\frac{x+y}{2},\tau)\vtb(\frac{x-y}{2},\tau)
\\
2\vtc(x,2\tau)\vtc(y,2\tau)=\vtc(\frac{x+y}{2},\tau)\vtc(\frac{x-y}{2},\tau)
+\vtd(\frac{x+y}{2},\tau)\vtd(\frac{x-y}{2},\tau)
\\
2\vtb(x,2\tau)\vtb(y,2\tau)=\vtc(\frac{x+y}{2},\tau)\vtc(\frac{x-y}{2},\tau)
-\vtd(\frac{x+y}{2},\tau)\vtd(\frac{x-y}{2},\tau)
\end{array}
\eq

\subsection{Appendix C. Deformed elliptic functions}
\setcounter{equation}{0}
\def\theequation{C.\arabic{equation}}

\beq{90}
\vf_a^\eta(z)=\bfe_N(a_2z)\phi(\frac{a_1+a_2\tau}{N}+\eta,z)\,,~
a\in\mZ^{(2)}_N\,,~\eta\in\Si_\tau\,.
\eq
It follows from (\ref{qp}) that $\vf_a^\eta(z)$ is well defined on
$\mZ^{(2)}_N$:
\beq{91}
\vf_{a+c}^\eta(z)=\vf_{a}^\eta(z)\,,~{\rm for}~c_{1,2}\in\mZ~{\rm mod}
\,N\,.
\eq
\beq{92}
\vf_a^\eta(z+1)=\bfe_N(a_2)\vf_a^\eta(z)\,,~~
\vf_a^\eta(z+\tau)=\bfe_N(-a_1-N\eta)\vf_a^\eta(z)\,.
\eq

The following formulae can be proved directly by checking the structure of
poles and quasi-periodic properties:

\beq{efi21}
\phi(w,\eta)\vfe_\al(z-w)+\phi(-w,\eta)\vfe_\al(z+w)=\phi(z,\eta)
(\vf_\al(z-w)+\vf_\al(z+w))
\eq

\beq{efi13}
\begin{array}{c}
\fe(z-w)\fe(z+w)+\vfe_\al(z-w)\vfe_\al(z+w)+\vfe_\be(z-w)\vfe_\be(z+w)
+\\
\vfe_\ga(z-w)\vfe_\ga(z+w)=
2\phi^{2\eta}(z-w)\phi^w(2\eta)+2\phi^{2\eta}(z+w)\phi^{-w}(2\eta)
\end{array}
\eq

\beq{efi14}
\begin{array}{c}
\fe(z-w)\fe(z+w)+\vfe_\al(z-w)\vfe_\al(z+w)-\vfe_\be(z-w)\vfe_\be(z+w)
-\\
\vfe_\ga(z-w)\vfe_\ga(z+w)=
2\phi^{2\eta}(z-w)\varphi_\al^w(2\eta)+
2\phi^{2\eta}(z+w)\varphi_\al^{-w}(2\eta)
\end{array}
\eq

\beq{efi15}
\begin{array}{c}
\vfe_\be(z-w)\fe(z+w)-\fe(z-w)\vfe_\be(z+w)-\vfe_\al(z-w)\vfe_\ga(z+w)
+\\
\vfe_\ga(z-w)\vfe_\al(z+w)=
2\varphi^{2\eta}_\be(z-w)\varphi^w_\al(2\eta)-
2\varphi^{2\eta}_\be(z+w)\varphi_\al^{-w}(2\eta)
\end{array}
\eq

\beq{efi16}
\begin{array}{c}
\vfe_\be(z-w)\vfe_\ga(z+w)+\vfe_\ga(z-w)\vfe_\be(z+w)-
\vfe_\al(z-w)\fe(z+w)
-\\
\fe(z-w)\vfe_\al(z+w)=
2\varphi^{2\eta}_\al(z-w)\varphi^w_\al(2\eta)+
2\varphi^{2\eta}_\al(z+w)\varphi_\al^{-w}(2\eta)
\end{array}
\eq

\beq{efi6}
\vfe_\be(z+w)\vfe_\gamma(z)\fe(w)+\vfe_\al(z+w)\fe(z)\vfe_\gamma(w)=
\eq
$$
\fe(z+w)\vfe_\al(z)\vfe_\be(w)+\vfe_\gamma(z+w)\vfe_\be(z)\vfe_\al(w)
$$

\beq{efi61}
\vfe_\be(z-w)\vfe_\gamma(z)\fe(w)-\vfe_\al(z-w)\fe(z)\vfe_\gamma(w)=
\eq
$$
-\fe(z-w)\vfe_\al(z)\vfe_\be(w)+\vfe_\gamma(z-w)\vfe_\be(z)\vfe_\al(w)
$$

\beq{efi7}
\begin{array}{c}
(E_1(\eta+\be)+E_1(\eta-\be)-E_1(\eta+\al)-E_1(\eta-\al))\times \\
(\vfe_\gamma(z+w)\vfe_\gamma(z)\fe(w)-\fe(z+w)\fe(z)\vfe_\gamma(w))
=\\
(E_1(\eta+\gamma)+E_1(\eta-\gamma)-2E_1(\eta))\times \\
(-\vfe_\al(z+w)\vfe_\al(z)\vfe_\be(w)+\vfe_\be(z+w)\vfe_\be(z)\vfe_\al(w))
\end{array}
\eq

\beq{efi8}
\begin{array}{c}
\fe(w)(-\vfe_\al(z-\om_\al)\vfe_\be(w-\om_\be)\vfe_\al(z-w)+
\vfe_\be(z-\om_\be)\vfe_\al(w-\om_\al)\vfe_\be(z-w))=
\\
-\phi^{-\eta}(w)(\vfe_\al(z-\om_\al)\vfe_\be(w-\om_\be)\vfe_\al(z+w)-
\vfe_\be(z-\om_\be)\vfe_\al(w-\om_\al)\vfe_\be(z+w))
\end{array}
\eq

\beq{efi9}
\begin{array}{c}
\fe(w+\om_\al)(\vfe_\ga(z-\om_\ga)\fe(w)\vfe_\be(z-w)-
\fe(z)\vfe_\ga(w-\om_\ga)\vfe_\al(z-w))=
\\
\fe(-w+\om_\al)(\vfe_\ga(z-\om_\ga)\fe(w)\vfe_\be(z+w)+
\fe(z)\vfe_\ga(w-\om_\ga)\vfe_\al(z+w))
\end{array}
\eq

\subsection{Appendix D. Comments to Proof of Reflection Equation}
\setcounter{equation}{0}
\def\theequation{D.\arabic{equation}}

Here we give some comments on the proof of the Proposition 3.2.

A direct substitution of (\ref{60}-\ref{61}) into (\ref{62})
yields three types expressions proportional to $\sigma\otimes 1$,
$1\otimes \sigma$ and $\sigma\otimes \sigma$. Consider, for example,
$1\otimes\sigma_\ga$ which contains the additional constants.
By the usage of (\ref{efi13}-\ref{efi16}) it simplifies to
\beq{w40}
\begin{array}{c}
[\Sh_\gamma,\Sh_0]\left(2\phi^w(\hbar)(
\vfh_\gamma(z)\fh(w)\vf_\gamma^{\hbar}(z-w)-
\fh(z)\vfh_\gamma(w)\phi^{\hbar}(z-w))\right.+\\
\left. 2\phi^{-w}(\hbar)(
\vfh_\gamma(z)\fh(w)\vf_\gamma^{\hbar}(z+w)-
\fh(z)\vfh_\gamma(w)\phi^{\hbar}(z+w))\right)+\\
\ [\tilde{\nu}_\gamma,\Sh_0]\left(2\phi^w(\hbar)(
\vfh_\gamma(z-\om_\ga)\fh(w)\vf_\gamma^{\hbar}(z-w)-
\fh(z)\vfh_\gamma(w-\om_{\gamma})\phi^{\hbar}(z-w))\right.+\\
\left. 2\phi^{-w}(\hbar)(
\vfh_\gamma(z-\om_\ga)\fh(w)\vf_\gamma^{\hbar}(z+w)-
\fh(z)\vfh_\gamma(w-\om_\ga)\phi^{\hbar}(z+w))\right)+\\
4i\tilde{\nu}_\al\Sh_\be\left(\vfh_\al(z-\om_\al)\vfh_\be(w)(
\vfh_\al(z-w)\phi(\hbar,w)+\vfh_\al(z+w)\phi(\hbar,-w))-\right.\\
\left.\vfh_\be(z)\vfh_\al(w-\om_\al)(
\vfh_\be(z-w)\phi(\hbar,w)+\vfh_\be(z+w)\phi(\hbar,-w))\right)+\\
4i\tilde{\nu}_\be\Sh_\al\left(\vfh_\al(z)\vfh_\be(w-\om_\be)(
\vfh_\al(z-w)\phi(\hbar,w)+\vfh_\al(z+w)\phi(\hbar,-w))-\right.\\
\left.\vfh_\be(z-\om_\be)\vfh_\al(w)(
\vfh_\be(z-w)\phi(\hbar,w)+\vfh_\be(z+w)\phi(\hbar,-w))\right)
=\\
i[\Sh_\al,\Sh_\be]_+\left(2\phi^w(\hbar)(
-\vfh_\al(z)\vfh_\be(w)\vf_\al^{\hbar}(z-w)+
\vfh_\be(z)\vfh_\al(w)\vf_\be^{\hbar}(z-w))\right.+\\
\left. 2\phi^{-w}(\hbar)(
-\vfh_\al(z)\vfh_\be(w)\vf_\al^{\hbar}(z+w)+
\vfh_\be(z)\vfh_\al(w)\vf_\be^{\hbar}(z+w))\right)+\\
\ i[\tilde{\nu}_\al,\tilde{\nu}_\be]_+\left(2\phi^w(\hbar)(
-\vfh_\al(z-\om_\al)\vfh_\be(w-\om_\be)\vf_\al^{\hbar}(z-w)+\right.\\
\left.\vfh_\be(z-\om_\be)\vfh_\al(w-\om_\al)\vf_\be^{\hbar}(z-w))\right.+\\
\left. 2\phi^{-w}(\hbar)(
-\vfh_\al(z-\om_\al)\vfh_\be(w-\om_\be)\vf_\al^{\hbar}(z+w)+
\vfh_\be(z-\om_\be)\vfh_\al(w-\om_\al)\vf_\be^{\hbar}(z+w))\right),
\end{array}
\eq
where $(\al,\be,\ga)$ is equivalent to $(1,2,3)$ under cyclic permutations.
The expression behind $[\tilde{\nu}_\al,\tilde{\nu}_\be]_+$ vanishes
due to (\ref{efi8}).
At the same time the expression behind $[\tilde{\nu}_\ga,\Sh_0]$ has a pole
at $w=\om_\ga$
different
from those of behind expressions  $[\Sh_\ga,\Sh_0]$ and
$[\Sh_\al,\Sh_\be]_+$. Thus $[\tilde{\nu}_\ga,\Sh_0]=0$.
At the moment we have
\beq{w401}
\begin{array}{c}
[\Sh_\gamma,\Sh_0]\left(2\vfh_\gamma(z)\fh(w)
(\phi^w(\hbar)\vf_\gamma^{\hbar}(z-w)+\phi^{-w}(\hbar)\vf^{\hbar}_\ga(z+w))
\right.-\\
\left.
2\fh(z)\vfh_\gamma(w)(\phi^w(\hbar)\phi^{\hbar}(z-w)+
\phi^{-w}(\hbar)\phi^{\hbar}(z+w))\right)\\
4i\tilde{\nu}_\al\Sh_\be\left(\vfh_\al(z-\om_\al)\vfh_\be(w)(
\vfh_\al(z-w)\phi(\hbar,w)+\vfh_\al(z+w)\phi(\hbar,-w))-\right.\\
\left.\vfh_\be(z)\vfh_\al(w-\om_\al)(
\vfh_\be(z-w)\phi(\hbar,w)+\vfh_\be(z+w)\phi(\hbar,-w))\right)+\\
4i\tilde{\nu}_\be\Sh_\al\left(\vfh_\al(z)\vfh_\be(w-\om_\be)(
\vfh_\al(z-w)\phi(\hbar,w)+\vfh_\al(z+w)\phi(\hbar,-w))-\right.\\
\left.\vfh_\be(z-\om_\be)\vfh_\al(w)(
\vfh_\be(z-w)\phi(\hbar,w)+\vfh_\be(z+w)\phi(\hbar,-w))\right)
=\\
i[\Sh_\al,\Sh_\be]_+\left(-2\vfh_\al(z)\vfh_\be(w)(
\phi^w(\hbar)\vf_\al^{\hbar}(z-w)+
\phi^{-w}(\hbar)\vf_\al^{\hbar}(z+w))\right.\\
\left.
2\vfh_\be(z)\vfh_\al(w)(
\phi^w(\hbar)\vf_\be^{\hbar}(z-w)+
\phi^{-w}(\hbar)\vf_\be^{\hbar}(z+w))\right)
\end{array}
\eq

Using then (\ref{efi21}) and cancelling $2\phi(z,\hbar)$ we have
\beq{w402}
\begin{array}{c}
[\Sh_\ga,\Sh_0]\left(\vfh_\ga(z)\fh(w)(\vf_\ga(z-w)+\vf_\ga(z+w))-\right.
\\
\left.\vfh_\ga(w)(\fh(w)\fh(z-w)+\fh(-w)\fh(z+w))\right)+
\\
2i\tilde{\nu}_\al\Sh_\be\left(\vfh_\al(z-\om_\al)\vfh_\be(w)
(\vf_\al(z-w)+\vf_\al(z+w))-\right.\\
\left.\vfh_\be(z)\vfh_\al(w-\om_\al)(\vf_\be(z-w)+\vf_\be(z+w))\right)+
\\
2i\tilde{\nu}_\be\Sh_\al\left(\vfh_\al(z)\vfh_\be(w-\om_\be)
(\vf_\al(z-w)+\vf_\al(z+w))-\right.\\
\left.\vfh_\be(z-\om_\be)\vfh_\al(w)(\vf_\be(z-w)+\vf_\be(z+w))\right)=
\end{array}
\eq
$$
=i[\Sh_\al,\Sh_\be]_+\left(-
\vfh_\al(z)\vfh_\be(w)(\vf_\al(z-w)+\vf_\al(z+w))+\right.\\
\left.\vfh_\be(z)\vfh_\al(w)(\vf_\be(z-w)+\vf_\be(z+w))\right)
$$
To get the final result one should compare the structure of poles
($w=0$ and $z=-w$).
(\ref{efi6}-\ref{efi9}).  Other types of expressions can be simplified
in the same way through the use of (\ref{efi6}-\ref{efi9}).

\small{

\end{document}
\begin{thebibliography}{60}

\bibitem{AL}\
 D.Arinkin, S.Lysenko,
 Isomorphisms between moduli spaces of SL(2)-bundles with connections
on ${\mathbb P}^1/\{x_1,\ldots,x_4\}$,
{\it Math. Res. Lett.}, {\bf 4},  (1997),   181--190, \\
 On the moduli spaces of SL(2)-bundles with connections
on ${\mathbb P}^1/\{x_1,\ldots,x_4\}$,
 {\it Internat. Math. Res. Notices}, no. 19, (1997),   983--999.

\bibitem{At}\
M.Atiyah,
Complex analytic connections in fibre bundles,
{\it Trans. of the Amer. Math. Soc.}, {\bf 85} (1957), 181--207.

\bibitem{Bax}\
R.J.Baxter,  Eight-vertex model in lattice statistics and one-dimensional
anisotropic Heisenberg chain, {\it I.Ann.Phys.}, {\bf 76}, (1973), 48--71.


\bibitem{BD}
A.~Belavin and V.~Drinfeld,  Solutions of
the classical Yang-Baxter equation for simple Lie algebras,
{\it Funct. Anal and Applic.}, {\bf 16},
  (1982), no. 3, 1--29.

\bibitem{BoMa}
\ A.V.Borisov, I.S.Mamaev,
{\em Modern Methods of the Theory of Integrable Systems},
Moscow - Izhevsk: Institute of Computer Science, (2003).

\bibitem{Ca}
F.Calogero,
Exactly solvable one dimensional many-body problem,
{\it Lettere al Nuovo Cimento}, {\bf 13}, (1975), 411-416.

\bibitem{FT2}\ L.Faddeev, L.Takhtajan, {\em Hamiltonian approach to solitons
theory}, Springer Series in Soviet Mathematics. Springer-Verlag, Berlin, (1987)

\bibitem{FO}
B.~Feigin and A.~Odesski,
 Sklyanin's elliptic algebras,
{\it Funct. Anal. and Applic.}, {\bf 23}, (1989), no. 3, 207--214.

\bibitem{Fedorov}\ Yu.N.Fedorov,
 Lax Representation with Spectral Parameter
on the Coverings of Hyperelliptic Curves,
{\it Math. notes},  {\bf 54}, (1993), 94--109.

\bibitem{Gam}\ B. Gambier,
 Sur les \'{e}quations differentielles du
 second ordre
et du premier degr\'{e} dont l'integral g\'{e}n\'{e}rale a ses points
critiques
fixes, {\it Acta Math. Ann}, {\bf 33}, (1910), 1-55.

\bibitem{GH}\
J.Gibbons and  T.Hermsen, A generalization of Calogero-Moser
system, {\it Physica D}, {\bf 11D}, (1984), 337--348.

\bibitem{Has}\ K.Hasegawa,  Ruijsenaars' commuting difference
operators as commuting transfer matrices,
{\it Comm. Math. Phys.},  {\bf 187},  (1997),  289--325.


\bibitem{H}\ N. Hitchin,  Stable  bundles and Integrable Systems,
{\it Duke Math. Journ.}, {\bf 54}, (1987), 91-114.

\bibitem{In}\
 V.Inozemtsev,
 Lax Representation with spectral
parameter on a torus for integrable particle systems,
{\it Lett. Math. Phys.},
{\bf 17}, (1989), 11-17.


\bibitem{InKo}\ T.Inami, H.Konno, Integrable XYZ spin chain with boundaries,
{\it J.Phys.A: Math.Gen.}, {\bf 27}, (1994), L913-L918.

\bibitem{Iw}\
 K.Iwasaki,
 Fuchsian moduli on a Riemann surface - its
Poisson structure and Poincar\'e-Lefschetz duality,
{\it Pacific J. Math.}, {\bf 155}, (1992), 319-340

\bibitem{KLO}\
B.Khesin, A.Levin, M.Olshanetsky,
Bihamiltonian structures and quadratic algebras in Hydrodynamics and on
non-commutative torus,
{\it Comm. Math. Phys.}, {\bf 270}, (2004), 581--612.

\bibitem{Kr}\
I.M.Krichever,
Elliptic solutions of Kadomtsev-Petviasvili equation and integrable system
of particles,
{\it Funct. Anal. Appl.} {\bf 14}, (1980), 282-290.

\bibitem{Kuznetsov}\
 V.B. Kuznetsov, M.F. Jorgensen, P.L. Christiansen,
New boundary conditions for integrable lattices,
{\it J.Phys. A} {\bf 28}, (1995), 4639-4654.


\bibitem{LO}\
 A.Levin, M.Olshanetsky,
 Hierarchies of isomonodromic deformations and Hitchin systems.
Moscow Seminar in Mathematical Physics,
{\it Amer. Math. Soc. Transl. Ser. 2}, {\bf 191},
 Amer. Math. Soc., Providence, RI, (1999), 223--262.

\bibitem{LO2}\
 A.Levin, M.Olshanetsky,
Non-autonomous Hamiltonian systems related to higher Hitchin integrals,
(Russian) {\it Teoret. Mat. Fiz. }, {\bf 123}, no. 2,  (2000),  237--263;
translation in  {\it Theoret. and Math. Phys.}, {\bf 123}  (2000), 609--632.

\bibitem{LOZ}\
 A.Levin, M.Olshanetsky, A.Zotov,
Hitchin systems--
symplectic Hecke correspondence and two-dimensional version,
{\it  Comm. Math. Phys.},  {\bf 236},  (2003),  93--133.

\bibitem{Ma}\
 Yu.I.Manin,
Sixth Painlev\'{e} equation, universal
elliptic curve, and mirror of $P^2$,
{\it Amer. Math. Soc. Transl. (2)}, {\bf 186}, (1998) 131-151.

\bibitem{Mo}\
  J.Moser,
 Three integrable Hamiltonian systems connected with isospectral
deformations,
{\it Advances in Math. } {\bf 16},  (1975), 197--220.

\bibitem{Mam}\
 D.Mumford, {\em Tata Lectures on Theta I, II,} Birkh\"{a}user
Boston, 1983, 1984.

\bibitem{N}\ N. Nekrasov, Holomorphic bundles and many-body
systems, {\it Comm. Math. Phys.}, {\bf 180} (1996) 587-604.

\bibitem{Ob}\
S.Oblezin,
Isomonodromic deformations and the Hecke correspondence,
{\it Moscow Math. Journ.}, {\bf 5}, (2005), 415-441.

\bibitem{O}\
M. Olshanetsky,  The large $N$ limits of integrable models, {\it Mosc. Math. J.},
 {\bf 3},  (2003),  1307--1331.

\bibitem{Pain}
P.Painlev\'{e},
Sur les \'{e}quations diff\'{e}rentielles du second odre \`{a}
points critics fixes,
{\it CRAS}, {\bf 143}, (1906), 1111-1117.


\bibitem{STSR}\
A.Reyman and M.Semenov-Tian-Schansky,
 Lie algebras and Lax equations with spectral parameter on elliptic curve,
(Russian) {\it Zap. Nauchn. Sem. Leningrad. Otdel. Mat. Inst. Steklov. (LOMI)}
{\bf 150}  (1986),  Voprosy Kvant. Teor. Polya i Statist. Fiz. 6, 104--118, 221;
translation in {\it J. Soviet Math.},  {\bf 46}, no. 1, (1989),   1631--1640.

\bibitem{Si}\
S.T.Simpson,  Harmonic bundles on non-compact curves,
{\it Journ. of AMS} {\bf 3} (1990), 713--770.

\bibitem{Skl2}\
E.Sklyanin,  Boundary conditions for integrable
equations, {\it Func. Anal. Appl.}, Vol.21, 86-87, (1987);

\bibitem{Skl}\
 E.Sklyanin, Some algebraic structures connected
with the Yang-Baxter equation,
{\it Funct. Anal. and Applic.}, {\bf 16} (1982), no. 4, 27--34.

\bibitem{Skl1}\
 E.Sklyanin, Boundary conditions for integrable quantum systems,
{\it J. Phys. A: Math. Gen}, {\bf 21}, (1988) 2375--2389.

\bibitem{SklTak1}
\ E.Sklyanin,  T.Takebe,  Algebraic Bethe ansatz for the $XYZ$ Gaudin model,
 {\it Phys. Lett. A}, {\bf 219}  (1996),  217--225.

\bibitem{SklTak}\
E.Sklyanin, T.Takebe,
  Separation of Variables in the Elliptic Gaudin Model,
 {\it Comm. Math. Phys.},  {\bf 204}  no. 1, (1999),  17--38.


\bibitem{Ta}\
 K.Takasaki,
 Spectral Curves and Whitham Equations in the
Isomonodromic Problems of Schlesinger Type,
{\it Asian J. Math.},  {\bf 2} no. 4, (1998),  1049--1078.


\bibitem{Va}\
V.Vakulenko,
Note on the Ruijsenaars-Schneider model,
  math.QA/9909079

\bibitem{Vol}\
 V.Volterra, Sur la th{\'e}orie des variations des
latitudes, {\it Acta Math.} {\bf 22} (1899),
201-357.

\bibitem{Wo}\
S.Wojciechowski,  An integrable marriage of the Euler equations
with the Calogero-Moser system, {\it Phys. Lett.}, {\bf A111} (1985) 101

\bibitem{LW}\
 Workshop on Painlev\'{e} Transcedents,
Their asymptotics and Physical Applications, {\it NATO ASI Ser. B:
Physics} Vol. 278 (1990) Sainte Adele,
Quebec, ed. by D. Levi and P. Winternitz


\bibitem{Zhuk}\
 N.E.Zhukovsky, {\it Jour. Phys. Chem. Soc.},
{\bf 17} (1885) 81-113, 145-199, 231-280.


\bibitem{Z}\
A.Zotov,
Elliptic Linear Problem for Calogero-Inozemtsev Model and Painlev{\'e}
VI Equation,
{\it Lett.Math.Phys.}, {\bf 67} (2004) 153-165.


\end{thebibliography}
